\documentclass[12pt]{article}
\usepackage{eurosym}
\usepackage{amsfonts}
\usepackage{amsmath,amssymb}
\usepackage[a4paper]{geometry}
\usepackage[utf8]{inputenc}
\usepackage[english]{babel}
\usepackage{amsthm}

\newtheorem{Lemma}{Lemma}[section]
\newtheorem{theorem}[Lemma]{Theorem}

\newtheorem{lemma}[Lemma]{Lemma}

\newtheorem{remark}[Lemma]{Remark}
\newtheorem{definition}[Lemma]{Definition}

\begin{document}
\begin{center}
\bigskip {\Large \textbf{Well-posedness of the Stochastic Degasperis-Procesi
Equation}} \vspace{3mm}\\[0pt]

\textsc{\ Lynnyngs K. Arruda}{\footnote{%
Departamento de Matem\'atica, Universidade Federal de S\~ao Carlos. CP 676,
13565-905, S\~ao Carlos - SP, Brazil. E-mail: lynnyngs{\char'100}dm.ufscar.br%
},}

\hspace{0mm} \textsc{\ Nikolai V. Chemetov}{\footnote{%
Department of Computing and Mathematics, University of S{\~a}o Paulo,
14040-901 Ribeir{\~a}o Preto - SP, Brazil, E-mail: nvchemetov{\char'100}%
gmail.com.}},

\hspace{0mm} \textsc{Fernanda Cipriano}{\footnote{%
Departamento de Matem\'atica, Faculdade de Ci\^encias e Tecnologia da
Universidade Nova de Lisboa and Centro de Matem\'atica e Aplica\c c\~oes,
Lisbon, Portugal E-mail: cipriano{\char'100}fct.unl.pt.}}
\end{center}

\vspace{3mm}
\date{\today }

\begin{abstract}
This article studies the Stochastic Degasperis-Procesi (SDP) equation on $%
\mathbb{R}$ with an additive noise. Applying the kinetic theory, and
considering the initial conditions in $L^2(\mathbb{R})\cap L^{2+\delta}(%
\mathbb{R})$, for arbitrary small $\delta>0$, we establish the existence of a
global pathwise solution. Restricting to the particular case of zero noise,
our result improves the deterministic solvability results that exist in the
literature.
\end{abstract}



\vspace{3mm}

\thispagestyle{empty}

\noindent 

\section{Introduction}

\setcounter{equation}{0} 
This work deals with the Degasperis-Procesi equation, which can be
considered as a particular case of the following ab-family of equations 
\begin{equation}
\partial _{t}u-\partial _{xxt}^{3}u+\partial _{x}(a(u,\partial
_{x}u))=\partial _{x}\left( b^{\prime }(u)\frac{(\partial _{x}u)^{2}}{2}%
+b(u)\partial _{xx}^{2}u\right) .  \label{Hyp1}
\end{equation}%
Namely, observing that 
\begin{equation*}
\partial _{x}\left( b(u)\partial _{x}u\right) =b^{\prime }(u)(\partial
_{x}u)^{2}+b(u)\partial _{xx}^{2}u,
\end{equation*}%
the right hand side of the equation (\ref{Hyp1}) corresponds to 
\begin{equation*}
\partial _{xx}^{2}\left( b(u)\partial _{x}u\right) -\partial _{x}\left(
b^{\prime }(u)\frac{(\partial _{x}u)^{2}}{2}\right) ,
\end{equation*}%
and the equation reads 
\begin{equation*}
\left( 1-\partial _{xx}^{2}\right) \partial _{t}u+\partial _{x}\left(
a(u,\partial _{x}u)-\int_{0}^{u}b(s)ds+b^{\prime }(u)\frac{(\partial
_{x}u)^{2}}{2}\right) +\left( 1-\partial _{xx}^{2}\right) \left(
b(u)\partial _{x}u\right) =0.
\end{equation*}%
Applying the operator $\left( 1-\partial _{xx}^{2}\right) ^{-1}$ to the
latter equation, we deduce the system 
\begin{equation}
\left\{ 
\begin{array}{l}
\partial _{t}u+\left( b(u)\partial _{x}u\right) +\partial _{x}p=0,\vspace{2mm%
} \\ 
\left( 1-\partial _{xx}^{2}\right) \mathrm{p}=a(u,\partial
_{x}u)-\int_{0}^{u}b(s)ds+b^{\prime }(u)\frac{(\partial _{x}u)^{2}}{2}.%
\vspace{2mm}%
\end{array}%
\right.   \label{Hyp2}
\end{equation}%
Considering the particular case 
\begin{equation}
b(u)=u,\qquad a(u,\partial _{x}u)=2u^{2}-\frac{(\partial _{x}u)^{2}}{2},
\label{choseab}
\end{equation}%
we recover the deterministic Degasperis-Procesi (DP) equation \cite{DP1998} 
\begin{equation}
\left\{ 
\begin{array}{l}
\partial _{t}u+u\partial _{x}u+\partial _{x}\mathrm{p}=0,\vspace{2mm} \\ 
\left( 1-\partial _{xx}^{2}\right) \mathrm{p}=\frac{3}{2}u^{2}.\vspace{2mm}%
\end{array}%
\right.   \label{Hyp22}
\end{equation}%
The DP equation was developed by the desire to find a water wave equation,
having traveling waves that break, in contrast to the Korteweg de Vries
equation which has only smooth solutions.

The DP equation \eqref{Hyp22} admits global weak solutions representing a
wave train of peakons \cite{LS1,LS2}. In addition it is completely
integrable, possesses a Lax pair, a bi-Hamiltonian structure, and an
infinite hierarchy of symmetries and conservation laws \cite{DHH}. The local
well-posedness of DP equation with initial data $u_{0}\in H^{s}(\mathbb{R})$%
, $s>3/2$, was proved in \cite{Yin1}. In \cite{Yin1}, it was also derived
the precise blow-up scenario and a blow-up result. In \cite{HHH} it was
shown that DP equation is ill-posed in Sobolev spaces $H^{s}(\mathbb{R})$
when $s<3/2.$ The global existence of strong solutions and global weak
solutions of DP equation were studied in \cite{LinLiu,Yin2,Yin3} and the
blow-up solutions in \cite{E1, E2, LiuYin, Yin3}.
 Global existence, $L^1$-stability and uniqueness results
for weak solutions in $L^{1}(\mathbb{R})\cap BV(\mathbb{R})$ and in $L^{2}(%
\mathbb{R})\cap L^{4}(\mathbb{R})$ with an additional entropy condition was
obtained in \cite{CK06} (see also \cite{Li}), and the periodic case was investigated in \cite{CK15}. 
The article  \cite{CK07} shows that the uniqueness result can be achieved by replacing Kruzkov's entropy conditions  by an Ole\v{ı}nik-type estimate.
 These weak solutions may
develop discontinuities in finite time. A detailed analysis in the articles 
\cite{E2}, \cite{K},\ \cite{L} explains the physical phenomenon of wave
breaking, arising from the collision of a peakon with an anti-peakon, where
the solution of the DP equations develop shock waves, having jump
discontinuities.

 The stochastic Degasperis-Procesi (SDP) equation has been considered for the first time in the article
 \cite{Chen015}, where the authors studied the global well-posedness of a stochastic dynamic, driven by a linear multiplicative noise, in the space 
 of sample paths $C([0,\infty),H^s(\mathbb{R})$, $s>3/2.$

 This article considers a stochastic dynamic driven by an additive noise, and is devoted to show the existence of a global  solution that  lives in a Lebesgue space 
of integrable functions, which can develops discontinuities in
finite time, being consistent with the physical interpretation of the wave breaking and shock-peakons phenomenon.  Within this very irregular framework,   it is not expected to have the pathwise uniqueness of the solution. Therefore, the usual methods to prove the existence of a stochastic strong solution  relying  on the pathwise uniqueness (see \cite{RT20} and references therein) can not be applied. 

We study the following model, where
the random fluctuations are modeled through an additive Wiener process
\begin{equation}
\left\{ 
\begin{array}{l}
du=-\left( u\partial _{x}u+\partial _{x}\mathrm{p}\right) \ dt+d{\mathcal{W}}%
_{t}, \\ 
\left( 1-\partial _{xx}^{2}\right) \mathrm{p}=\frac{3}{2}u^{2}.\vspace{2mm}%
\end{array}%
\right.   \label{Hyp2S}
\end{equation}%
We establish the existence of a global 
 stochastic entropy solution for  \eqref{Hyp2S}
  with sample paths in $L^{\infty }(0,T;L^{2}(\mathbb{R}))$, in the sense of Definition \ref{DEF}.
  
 In order to capture  the relevant
physical property of the SDP equation corresponding to shock-peakon, we seek solutions with less regularity as possible, by considering the initial conditions just in $L^{2}(\mathbb{R})\cap L^{2+\delta }(\mathbb{R})$ for arbitrary small $\delta >0$. Our method is based on a conjugation of the deterministic entropy method,
 with a stochastic pathwise approach, adopting the strategy introduced 
by \cite{Fla99} to study the stochastic Euler equation. More precisely, we show the existence of a pathwise entropy solution, which corresponds to a stochastic process with values in the  $L^{\infty }(0,T;L^{2}(\mathbb{R})\cap L^{2+\delta
}(\mathbb{R}))$. The proof of the pathwise uniqueness for initial conditions in $L^{2}(\mathbb{R})\cap L^{2+\delta }(\mathbb{R})$ is not an easy task, and remains an open problem.  

We shoul mention that our solutions are strong in the
probability sense, i.e. for a Wiener process $\mathcal{W}_{t}$ given in advance, not being part of the solution. 
 We hope that the pathwise solution constructed in the present work will be relevant to study the statistical properties of the shock-peakons phenomenon, in  Degasteris-Process dynamic.
 
As far as we know, this is the first time the existence of the solution for the SDP is established by assuming only integrability for the initial conditions. Furthermore, as a particular case, $\mathcal{W}_{t}$ can be considered as a deterministic external force which is a distributional derivative of a time-continuous function, and our result in new even in the deterministic setting, since it extends the existence deterministic result in \cite{CK06}, where the existence is shown for  the initial conditions  in  $L^{2}(\mathbb{R})\cap L^{4}(\mathbb{R})$.

This paper is organized as follows. In Section 2, we formulate the problem. In
Section 3 we provide the proofs of technical results on the regularized
problem and a priori estimates. Finally, Section 4 is concerned with the
limit transition by the kinetic method.

\bigskip

\section{Formulation of the problem}

\setcounter{equation}{0}

Let us define some normed spaces. We denote the norm in $L^{s}(\mathbb{R})$
for $s\geqslant 1$\ by the associated norm by $\Vert \cdot \Vert _{L^{s}}.$
Let $X$ be a real Banach space with norm $\left\Vert \cdot \right\Vert _{X}.$
We denote by $L^{r}(0,T;X)$ the space of $X$-valued measurable $r-$%
integrable functions defined on $[0,T]$ for $r\geqslant 1$.

Let $(\Omega, \mathcal{F}, P)$ be a probability space. For $q>2$, we
consider a stochastic process $\mathcal{W}_t$, $t\geq 0$, defined on $(\Omega, \mathcal{F}, P)$, with values in the Banach space $L^{2}(%
\mathbb{R})\cap L^{q}(\mathbb{R})$ (for instance an infinite dimensional
Wiener process, see \cite{DaPrato}) with the following additional regularity
for $P-$a.e. $\omega\in\Omega$ 
\begin{equation*}
\partial_x\mathcal{W}(\cdot,\omega,\cdot)\in L^2(0,T;L^\infty(\mathbb{R})).
\end{equation*}
We define the filtation $\{\mathcal{F}_t\}_{t\geq 0}$, with 
$\mathcal{F}_t=\sigma\left(\mathcal{W}_s:\,0\leq s\leq t\right).$

\bigskip This work is devoted to the study of the well-posedness of the
stochastic Degasperis-Procesi equation in the one-dimensional domain $%
\mathbb{R}$, with an additive noise. The evolution equations are given by%
\begin{equation}
\left\{ 
\begin{array}{l}
du=-\left( u\partial _{x}u+\partial _{x}\mathrm{p}\right) \ dt+d{\mathcal{W}}%
_{t},\quad \left( 1-\partial _{xx}^{2}\right) \mathrm{p}=\frac{3}{2}%
u^{2}\quad \text{in }\mathbb{R}_{T}=(0,T)\times \mathbb{R}, \\ 
\\ 
u(0)=u_{0}\qquad \qquad \qquad \text{ in }\mathbb{R}%
\end{array}%
\right.   \label{Hyp3}
\end{equation}%
where $u=u(t,x)$ denotes the velocity of the fluid and $u_{0}=u_{0}(x)$ is
an initial condition.

The solution of this problem is understood in the following sense.

\begin{definition}
\label{DEF} A stochastic process $(u,\mathrm{p})$\ is a stochastic entropy
solution of (\ref{Hyp3}) if for $P-$a.e. $\omega \in \Omega $ \ this pair%
\begin{equation*}
u(\cdot ,\omega )\in L^{\infty }(0,T;L^{2}(\mathbb{R})),\qquad ~\mathrm{p}%
(\cdot ,\omega )\in L^{\infty }(0,T;W^{2,1}(\mathbb{R})),
\end{equation*}%
satisfy the equation 
\begin{equation}
\left( 1-\partial _{xx}^{2}\right) \mathrm{p}=\frac{3}{2}u^{2}\qquad \text{
a.e. in }\mathbb{R}_{T},  \label{1O}
\end{equation}%
and $z(\cdot ,\omega )=u(\cdot ,\omega )-\mathcal{W}(\cdot ,\omega )$ and $%
\mathrm{p}(\cdot ,\omega )$ fulfill 
\begin{equation}
\begin{array}{l}
\displaystyle\int_{\mathbb{R}_{T}}\biggl\{|z-c|\partial _{t}\varphi +\mathrm{%
sign}(z-c)\left( \left( \frac{z^{2}}{2}-\frac{c^{2}}{2}\right) +\left(
z-c\right) \mathcal{W}\right) \partial _{x}\varphi \vspace{1mm} \\ 
-\left( \partial _{x}\mathrm{p}+\mathcal{W}\partial _{x}\mathcal{W}%
+c\partial _{x}\mathcal{W}\right) \mathrm{sign}(z-c)\varphi \vspace{1mm}%
\biggr\}dtdx\vspace{1mm} \\ 
+\displaystyle\int_{\mathbb{R}}|z_{0}-c|\varphi (0,\text{\textbf{%
\textperiodcentered }})\,\,dx\geqslant 0%
\end{array}
\label{pdef}
\end{equation}%
for arbitrary $c\in \mathbb{R}$ \ and for any positive function $\varphi
(t,x)\in C^{\infty }([0,T];C_{0}^{\infty }(\mathbb{R}))$, such that $\varphi
(T,\text{\textperiodcentered })=0$, where $\ z_{0}=u_{0}.$
\end{definition}

\bigskip

Next, we state the main result of the article.

\begin{theorem}
\label{the_1} Let $\mathcal{W}$ be a stochastic process, such that for $P-$%
a.e. $\omega \in \Omega $ 
\begin{eqnarray}
\mathcal{W}(\cdot ,\omega ) &\in &C([0,T];L^{2}(\mathbb{R})\cap L^{q}(%
\mathbb{R})),\qquad \partial _{x}\mathcal{W}(\cdot ,\omega )\in
L^{2}(0,T;L^{\infty }(\mathbb{R})),  \notag \\
u_{0} &\in &L^{2}(\mathbb{R})\cap L^{q}(\mathbb{R})\qquad \text{with }q>2%
\text{.}  \label{au}
\end{eqnarray}%
Then there exists at least one stochastic entropy solution $(u,\mathrm{p})$\
of system (\ref{Hyp3}), such that for $P-$a.e. $\omega \in \Omega $ 
\begin{eqnarray*}
u(\cdot ,\omega ) &\in &L^{\infty }(0,T;L^{2}(\mathbb{R})\cap L^{q}(\mathbb{R%
})), \\
\mathrm{p}(\cdot ,\omega ) &\in &L^{\infty }(0,T;W^{1,s}(\mathbb{R})),\qquad
\forall s\in \left[ 1,+\infty ,\right]  \\
\partial _{xx}^{2}\mathrm{p}(\cdot ,\omega ) &\in &L^{\infty }(0,T;L^{q/2}(%
\mathbb{R}))\qquad \text{and}\qquad \mathrm{p}\geqslant 0\quad \text{a.e. in
\ }\mathbb{R}_{T}.
\end{eqnarray*}
\end{theorem}

\bigskip

\bigskip

\section{Existence of the regularized problem and a priori estimates}

\bigskip

\setcounter{equation}{0}

Let us introduce $z=u-\mathcal{W}$. By (\ref{Hyp3}) for $P-$a.e. $\omega \in
\Omega $, the stochastic processes $z$ and $\mathrm{p}$ should satisfy the
Cauchy problem 
\begin{equation*}
\left\{ 
\begin{array}{l}
\partial _{t}z=-(z+\mathcal{W})\partial _{x}(z+\mathcal{W})-\partial _{x}%
\mathrm{p}\qquad \text{in }\mathbb{R}_{T}, \\ 
\\ 
\left( 1-\partial _{xx}^{2}\right) \mathrm{p}=\frac{3}{2}(z+\mathcal{W}%
)^{2}\qquad \text{in }\mathbb{R}_{T}, \\ 
\\ 
z(0)=z_{0}=u_{0}\qquad \qquad \qquad \text{ in }\mathbb{R}.%
\end{array}%
\right. 
\end{equation*}%
In order to show the existence and uniqueness of the solution for the above
system, we follow the viscosity approach by introducing for any $\varepsilon
\in (0,1)$ \ the following stochastic viscosity system 
\begin{equation}
\left\{ 
\begin{array}{l}
\partial _{t}z_{\varepsilon }=-(z_{\varepsilon }+\mathcal{W})\partial
_{x}(z_{\varepsilon }+\mathcal{W})-\partial _{x}\mathrm{p}_{\varepsilon
}+\varepsilon \partial _{xx}^{2}z_{\varepsilon }\qquad \text{in }\mathbb{R}%
_{T}, \\ 
\\ 
\left( 1-\partial _{xx}^{2}\right) \mathrm{p}_{\varepsilon }=\frac{3}{2}%
(z_{\varepsilon }+\mathcal{W})^{2}\qquad \text{in }\mathbb{R}_{T}, \\ 
\\ 
z_{\varepsilon }(0)=z_{\varepsilon ,0}=u_{\varepsilon ,0}\qquad \qquad
\qquad \text{ in }\mathbb{R}.%
\end{array}%
\right.   \label{Sys}
\end{equation}%
Having that $u_{0}\in L^{2}(\mathbb{R})\cap L^{q}(\mathbb{R})$ with $q$
defined in (\ref{au}),\ using the regularization procedure we can assume $%
u_{\varepsilon ,0}\in C_{0}^{\infty }(\mathbb{R}),$ such that 
\begin{eqnarray*}
\Vert u_{\varepsilon ,0}\Vert _{L^{2}(\mathbb{R})} &\mathbf{\leqslant }%
&C\Vert u_{0}\Vert _{L^{2}(\mathbb{R})},\qquad \Vert u_{\varepsilon ,0}\Vert
_{L^{q}(\mathbb{R})}\mathbf{\leqslant }C\Vert u_{0}\Vert _{L^{q}(\mathbb{R}%
)}, \\
u_{\varepsilon ,0} &\rightarrow &u_{\varepsilon ,0}\quad \text{strongly in }%
\quad L^{2}(\mathbb{R})\cap L^{q}(\mathbb{R}).
\end{eqnarray*}%
For $P-$a.e. $\omega \in \Omega $,\ \ the existence of the unique solution $%
z_{\varepsilon }(\cdot ,\omega )\in H^{1}(0,T;L^{2}(\mathbb{R}))\cap
L^{2}(0,T;H^{2}(\mathbb{R}))$ \ for the system (\ref{Sys}) can be show by
standard results for a coupled parabolic - elliptic systems, for details we
refer to \cite{CK}, \cite{Evans}, \cite{LS}.

\subsection{$\protect\varepsilon -$Uniform estimates}

\bigskip

In this section we establish some $\varepsilon -$uniform estimates for the
pair $(z_{\varepsilon },\mathrm{p}_{\varepsilon }).$ First for $P-$a.e. $%
\omega \in \Omega $ let us define the solution $\ y_{\varepsilon }(\cdot
,\omega )$ of the equation 
\begin{equation}
(4-\partial _{xx}^{2})y_{\varepsilon }=z_{\varepsilon }\qquad \text{in }%
\mathbb{R}_{T}  \label{y}
\end{equation}%
and consider that $\mathrm{p}_{\varepsilon }=\mathrm{p}_{\varepsilon }^{(1)}+%
\mathrm{p}_{\varepsilon }^{(2)}$ with $\mathrm{p}_{\varepsilon }^{(1)}(\cdot
,\omega ),$ $\mathrm{p}_{\varepsilon }^{(2)}(\cdot ,\omega )$ being the
solutions \ of the equations 
\begin{equation}
\left( 1-\partial _{xx}^{2}\right) \mathrm{p}_{\varepsilon }^{(1)}=\frac{3}{2%
}z_{\varepsilon }^{2},\qquad \left( 1-\partial _{xx}^{2}\right) \mathrm{p}%
_{\varepsilon }^{(2)}=\frac{3}{2}\left( 2z_{\varepsilon }\mathcal{W}+%
\mathcal{W}^{2}\right) \qquad \text{in }\mathbb{R}_{T}.  \label{p}
\end{equation}%
Also we introduce the quantity 
\begin{equation*}
S(u)=4\Vert u\Vert _{L^{2}}^{2}+5\Vert \partial _{x}u\Vert
_{L^{2}}^{2}+\Vert \partial _{xx}^{2}u\Vert _{L^{2}}^{2}.
\end{equation*}

The following lemma collects the $\varepsilon -$uniform estimates for $%
(z_{\varepsilon },\mathrm{p}_{\varepsilon }).$

\begin{lemma}
\label{lem1} Under conditions (\ref{au}), for $P-$a.e. $\omega \in \Omega $
\ there exist some constants $C=C(\omega ),$ which are independent of $%
\varepsilon $, but may depend of $\omega ,$ such that%
\begin{eqnarray}
\Vert \sqrt{\varepsilon }\partial _{x}z_{\varepsilon }\Vert _{L^{2}(\mathbb{R%
}_{T})}^{2} &<&C,  \notag \\
\Vert z_{\varepsilon }\Vert _{L^{\infty }(0,T;L^{2}(\mathbb{R}))}
&<&C,\qquad \Vert z_{\varepsilon }\Vert _{L^{\infty }(0,T;L^{q}(\mathbb{R}))}%
\mathbf{<}C,\text{ }  \notag \\
\Vert \mathrm{p}_{\varepsilon }\Vert _{L^{\infty }(0,T;L^{p}(\mathbb{R}%
))},~\Vert \partial _{x}\mathrm{p}_{\varepsilon }\Vert _{L^{\infty
}(0,T;L^{p}(\mathbb{R}))} &<&C,\qquad \forall p\in \left[ 1,+\infty \right] ,
\notag \\
\Vert \partial _{xx}^{2}\mathrm{p}_{\varepsilon }\Vert _{L^{\infty
}(0,T;L^{q/2}(\mathbb{R}))} &<&C  \label{Z}
\end{eqnarray}%
with $q$ defined in (\ref{au}).
\end{lemma}

\textbf{Proof. \ }The equation (\ref{Sys})$_{1}$ \ in terms of $\mathrm{p}%
_{\varepsilon }^{(1)},\mathrm{p}_{\varepsilon }^{(2)}$ and $z_{\varepsilon }$
\ is written as the equality%
\begin{eqnarray}
\partial _{t}z_{\varepsilon }-\varepsilon \partial _{xx}^{2}z_{\varepsilon }
&=&-z_{\varepsilon }\partial _{x}z_{\varepsilon }-\partial _{x}\mathrm{p}%
_{\varepsilon }^{(1)}  \notag \\
&&-z_{\varepsilon }\partial _{x}\mathcal{W}-\mathcal{W}\partial
_{x}z_{\varepsilon }-\mathcal{W}\partial _{x}\mathcal{W}-\partial _{x}%
\mathrm{p}_{_{\varepsilon }}^{(2)}.  \label{zz}
\end{eqnarray}%
First note that 
\begin{equation*}
\int_{\mathbb{R}}\left( \partial _{t}z_{\varepsilon }-\varepsilon \partial
_{xx}^{2}z_{\varepsilon }\right) \left( y_{\varepsilon }-\partial
_{xx}^{2}y_{\varepsilon }\right) \ dx=\frac{1}{2}\frac{d}{dt}%
S(y_{\varepsilon })+\varepsilon S(\partial _{x}y_{\varepsilon })
\end{equation*}%
and%
\begin{equation*}
\int_{\mathbb{R}}\left( z_{\varepsilon }\partial _{x}z_{\varepsilon
}+\partial _{x}\mathrm{p}_{\varepsilon }^{(1)}\right) (y_{\varepsilon
}-\partial _{xx}^{2}y_{\varepsilon })\ dx=0.
\end{equation*}%
Hence, the multiplication of equality (\ref{zz}) by $\left( y_{\varepsilon
}-\partial _{xx}^{2}y_{\varepsilon }\right) $ yields 
\begin{align}
& \frac{d}{dt}S(y_{\varepsilon })+2\varepsilon S(\partial _{x}y_{\varepsilon
})  \label{a1} \\
& =2\int_{\mathbb{R}}\left( -z_{\varepsilon }\partial _{x}\mathcal{W}-%
\mathcal{W}\partial _{x}z_{\varepsilon }-\mathcal{W}\partial _{x}\mathcal{W}%
-\partial _{x}\mathrm{p}_{\varepsilon }^{(2)}\right) \left( y_{\varepsilon
}-\partial _{xx}^{2}y_{\varepsilon }\right) \ dx.  \notag
\end{align}%
Now, we estimate the right hand side in this identity. By (\ref{y}) we have%
\begin{equation}
\biggl|\int_{\mathbb{R}}z_{\varepsilon }\partial _{x}\mathcal{W}\left(
y_{\varepsilon }-\partial _{xx}^{2}y_{\varepsilon }\right) \ dx\biggr|%
\mathbf{\leqslant }C\Vert \partial _{x}\mathcal{W}\Vert _{L^{\infty
}}S(y_{\varepsilon }).  \label{a2}
\end{equation}%
Integrating by parts, we derive 
\begin{align*}
& \int_{\mathbb{R}}\left( -\mathcal{W}\partial _{x}z_{\varepsilon }\right)
\left( y_{\varepsilon }-\partial _{xx}^{2}y_{\varepsilon }\right) \ dx \\
& =2\int_{\mathbb{R}}\partial _{x}\mathcal{W}(y_{\varepsilon })^{2}\
dx-2\int_{\mathbb{R}}\partial _{x}\mathcal{W}(\partial _{x}y_{\varepsilon
})^{2}\ dx-\int_{\mathbb{R}}\partial _{xx}^{2}y_{\varepsilon }\partial _{x}(%
\mathcal{W}y_{\varepsilon })\ dx \\
& \qquad +\frac{1}{2}\int_{\mathbb{R}}\partial _{x}\mathcal{W}\left(
\partial _{xx}^{2}y_{\varepsilon }\right) ^{2}\ dx \\
& =2\int_{\mathbb{R}}\partial _{x}\mathcal{W}\ (y_{\varepsilon })^{2}\ dx-%
\frac{3}{2}\int_{\mathbb{R}}\partial _{x}\mathcal{W}\ (\partial
_{x}y_{\varepsilon })^{2}\ dx-\int_{\mathbb{R}}\partial
_{xx}^{2}y_{\varepsilon }\ \partial _{x}\mathcal{W}y_{\varepsilon }\ dx \\
& \qquad +\frac{1}{2}\int_{\mathbb{R}}\partial _{x}\mathcal{W}\ \left(
\partial _{xx}^{2}y_{\varepsilon }\right) ^{2}\ dx
\end{align*}%
that gives 
\begin{equation}
\biggl|\int_{\mathbb{R}}\left( -\mathcal{W}\partial _{x}z_{\varepsilon
}\right) \left( y_{\varepsilon }-\partial _{xx}^{2}y_{\varepsilon }\right) \
dx\biggr|\mathbf{\leqslant }C\Vert \partial _{x}\mathcal{W}\Vert _{L^{\infty
}}S(y_{\varepsilon }).  \label{a3}
\end{equation}%
In addition, we have the following estimates 
\begin{align}
\biggl|\int_{\mathbb{R}}\left( \mathcal{W}\partial _{x}\mathcal{W}\right)
\left( y_{\varepsilon }-\partial _{xx}^{2}y_{\varepsilon }\right) \ dx\biggr|%
& \mathbf{\leqslant }\Vert \mathcal{W}\Vert _{L^{2}}\Vert \partial _{x}%
\mathcal{W}\Vert _{L^{\infty }}\Vert y_{\varepsilon }-\partial
_{xx}^{2}y_{\varepsilon }\Vert _{L^{2}}  \notag \\
& \mathbf{\leqslant }\Vert \mathcal{W}\Vert _{L^{2}}^{2}\Vert \partial _{x}%
\mathcal{W}\Vert _{L^{\infty }}^{2}+CS(y_{\varepsilon }),  \label{a33}
\end{align}%
and%
\begin{align}
& \biggl|\int_{\mathbb{R}}\left( \partial _{x}\mathrm{p}_{\varepsilon
}^{(2)}\right) \left( y_{\varepsilon }-\partial _{xx}^{2}y_{\varepsilon
}\right) \ dx\biggr|=\biggl|\int_{\mathbb{R}}\left( \mathrm{p}_{\varepsilon
}^{(2)}-\partial _{xx}^{2}\mathrm{p}_{\varepsilon }^{(2)}\right) \partial
_{x}y_{\varepsilon }\ dx\biggr|  \notag \\
& =\frac{3}{2}\biggl|\int_{\mathbb{R}}\left( 2(4y_{\varepsilon }-\partial
_{xx}^{2}y_{\varepsilon })\mathcal{W}+\mathcal{W}^{2}\right) \partial
_{x}y_{\varepsilon }\ dx\biggr|  \notag \\
& =\frac{3}{2}\biggl|\int_{\mathbb{R}}(-4(y_{\varepsilon })^{2}+(\partial
_{x}y_{\varepsilon })^{2})\ \partial _{x}\mathcal{W}-2\mathcal{W}\partial
_{x}\mathcal{W}y_{\varepsilon }\ dx\biggr|  \notag \\
& \mathbf{\leqslant }C\Vert \partial _{x}\mathcal{W}\Vert _{L^{\infty
}}S(y_{\varepsilon })+\Vert \partial _{x}\mathcal{W}\Vert _{L^{\infty
}}\Vert \mathcal{W}\Vert _{L^{2}}\Vert y_{\varepsilon }\Vert _{L^{2}}  \notag
\\
& \mathbf{\leqslant }C\left( \Vert \partial _{x}\mathcal{W}\Vert _{L^{\infty
}}+1\right) S(y_{\varepsilon })+\Vert \partial _{x}\mathcal{W}\Vert
_{L^{\infty }}^{2}\Vert \mathcal{W}\Vert _{L^{2}}^{2}.  \label{a333}
\end{align}%
Combining \eqref{a1}-\eqref{a333}, we deduce 
\begin{equation*}
\frac{d}{dt}S(y_{\varepsilon })+2\varepsilon S(\partial _{x}y_{\varepsilon })%
\mathbf{\leqslant }\alpha (t)S(y_{\varepsilon })+\beta (t)
\end{equation*}%
with 
\begin{align*}
& \alpha (t)=C\left( \Vert \partial _{x}\mathcal{W}\Vert _{L^{\infty
}}+1\right) \in L^{1}(0,T), \\
& \beta (t)=C\Vert \mathcal{W}\Vert _{L^{2}}^{2}\Vert \partial _{x}\mathcal{W%
}\Vert _{L^{\infty }}^{2}\in L^{1}(0,T)\qquad \text{a.s. }\omega \in \Omega .
\end{align*}%
Applying the Gronwall's Lemma, for a.s. $\omega \in \Omega $ we obtain the
path-wise inequality 
\begin{equation}
S(y_{\varepsilon }(t))+2\varepsilon \int_{0}^{t}S(\partial
_{x}y_{\varepsilon })\ ds\mathbf{\leqslant }\left( S(y_{\varepsilon
}(0))+A\right) B\qquad \text{a.e. }t\in \lbrack 0,T]  \label{ss}
\end{equation}%
by the assumptions \eqref{au}. Here we denote $A(\omega )=\int_{0}^{T}\beta
(s)\,ds,\qquad B(\omega )=\exp \left( \int_{0}^{T}\alpha (s)\,ds\right) $.
This estimate gives \eqref{Z}$_{1}.$

Since $y_{\varepsilon },$ $z_{\varepsilon }$ fulfill \eqref{y}, there exist
positive constants $C_{1},C_{2},$ such that 
\begin{equation*}
\Vert z_{\varepsilon }(t)\Vert _{L^{2}}^{2}\leqslant 2\Vert y_{\varepsilon
}(t)\Vert _{L^{2}}^{2}+32\Vert \partial _{xx}^{2}y_{\varepsilon }(t)\Vert
_{L^{2}}^{2}\leqslant C_{1}S(y_{\varepsilon }(t)),\qquad \forall t\in
\lbrack 0,T]
\end{equation*}%
and the integration by parts gives%
\begin{equation*}
\Vert z_{\varepsilon }(t)\Vert _{L^{2}}^{2}=\Vert y_{\varepsilon }(t)\Vert
_{L^{2}}^{2}+8\Vert \partial _{x}y_{\varepsilon }(t)\Vert
_{L^{2}}^{2}+16\Vert \partial _{xx}^{2}y_{\varepsilon }(t)\Vert
_{L^{2}}^{2}\geqslant C_{2}S(y_{\varepsilon }(t)).
\end{equation*}%
These inequalities and \eqref{ss} imply 
\begin{eqnarray*}
\Vert z_{\varepsilon }(t)\Vert _{L^{2}}^{2} &\leqslant &CS(y_{\varepsilon
}(t))\leqslant C\left( S(y_{\varepsilon }(0))+1\right) \mathbf{\leqslant }%
C(\Vert z_{\varepsilon }(0)\Vert _{L^{2}}^{2}+1) \\
&\mathbf{\leqslant }&C(\Vert u_{\varepsilon ,0}\Vert _{L^{2}}^{2}+\Vert 
\mathcal{W}\Vert _{L^{2}}^{2}+1)\mathbf{\leqslant }C(\omega )\qquad \text{%
a.s. }\omega \in \Omega ,
\end{eqnarray*}%
which is \eqref{Z}$_{2}$.

Let us show the estimates \eqref{Z}$_{4,5}$ for the function $\mathrm{p}%
_{\varepsilon }.$ Taking into account that $\mathrm{p}_{\varepsilon }=%
\mathrm{p}_{\varepsilon }^{(1)}+\mathrm{p}_{\varepsilon }^{(2)}$ has the
following explicit representation%
\begin{equation}
\mathrm{p}_{\varepsilon }(t,x)=\frac{3}{4}\int_{\mathbb{R}}(z_{\varepsilon
}(t,y)+\mathcal{W}(t,y))^{2}\exp (-|x-y|)\ dy,\qquad (t,x)\in \mathbb{R}_{T}
\label{PP}
\end{equation}%
as the solution of \eqref{Sys}$_{2}$, we derive that%
\begin{equation}
\partial _{x}\mathrm{p}_{\varepsilon }(t,x)=-\frac{3}{4}\int_{\mathbb{R}%
}(z_{\varepsilon }(t,y)+\mathcal{W}(t,y))^{2}\exp (-|x-y|)sign(x-y)\ dy
\label{PP1}
\end{equation}%
and%
\begin{eqnarray*}
\Vert \mathrm{p}_{\varepsilon }\Vert _{L^{\infty }},~\Vert \partial _{x}%
\mathrm{p}_{\varepsilon }\Vert _{L^{\infty }} &\mathbf{\leqslant }&C\int_{%
\mathbb{R}}(z_{\varepsilon }(t,y)+\mathcal{W}(t,y))^{2}\ dy \\
&\leqslant &C\left( \Vert z_{\varepsilon }\Vert _{L^{2}}^{2}+\Vert \mathcal{W%
}\Vert _{L^{2}}^{2}\right) \mathbf{\leqslant }C(\omega ),\qquad \forall t\in
\lbrack 0,T].
\end{eqnarray*}%
Moreover Fubini's theorem gives 
\begin{eqnarray*}
\Vert \mathrm{p}_{\varepsilon }\Vert _{L^{1}},\Vert \partial _{x}\mathrm{p}%
_{\varepsilon }\Vert _{L^{1}} &\leqslant &C\int_{\mathbb{R}}\left( \int_{%
\mathbb{R}}(z_{\varepsilon }(t,y)+\mathcal{W}(t,y))^{2}\exp (-|x-y|)\
dx\right) dy \\
&\leqslant &C\left( \Vert z_{\varepsilon }\Vert _{L^{2}}^{2}+\Vert \mathcal{W%
}\Vert _{L^{2}}^{2}\right) \mathbf{\leqslant }C(\omega ),\qquad \forall t\in
\lbrack 0,T].
\end{eqnarray*}%
Applying the deduced estimates in $L^{\infty }$ and $L^{1},$ for any $p\in
\left( 1,+\infty \right) $\ we have 
\begin{eqnarray*}
\Vert \mathrm{p}_{\varepsilon }\Vert _{L^{p}}^{p} &\leqslant &\Vert \mathrm{p%
}_{\varepsilon }\Vert _{L^{p}}^{p-1}\Vert \mathrm{p}_{\varepsilon }\Vert
_{L^{1}}\mathbf{\leqslant }C(\omega ,p), \\
\Vert \partial _{x}\mathrm{p}_{\varepsilon }\Vert _{L^{p}}^{p} &\leqslant
&\Vert \partial _{x}\mathrm{p}_{\varepsilon }\Vert _{L^{p}}^{p-1}\Vert
\partial _{x}\mathrm{p}_{\varepsilon }\Vert _{L^{1}}\mathbf{\leqslant }%
C(\omega ,p),\qquad \forall t\in \lbrack 0,T],
\end{eqnarray*}%
that gives \eqref{Z}$_{4,5}.$

Now let us show \eqref{Z}$_{3}.$\ Now defining the convex function $\eta
(z)=\left( z^{2}\right) ^{q/2}$, we have $\eta (z_{\varepsilon })\equiv
z_{\varepsilon }{}^{q}=|z_{\varepsilon }|^{q}$ and $\eta ^{\prime
}(z_{\varepsilon })=qz_{\varepsilon }{}^{q-1}\equiv q|z_{\varepsilon }|^{q-1}%
\mathrm{sign}z_{\varepsilon }.$ Multiplying the equation (\ref{Sys})$_{1}$
by $\eta ^{\prime }(z_{\varepsilon })$, we deduce 
\begin{align*}
& \partial _{t}\left( |z_{\varepsilon }|^{q}\right) +\frac{q}{q+1}\partial
_{x}\left( z_{\varepsilon }^{q+1}\right) +q\partial _{x}\mathrm{p}%
_{\varepsilon }|z_{\varepsilon }|^{q-1}\mathrm{sign}z_{\varepsilon } \\
& =\varepsilon \partial _{x}\left( \partial _{x}|z_{\varepsilon
}|^{q}\right) -\varepsilon (\partial _{x}z_{\varepsilon })^{2}\eta ^{\prime
\prime }(z_{\varepsilon }) \\
& -q\partial _{x}\mathcal{W}|z_{\varepsilon }|^{q}-\mathcal{W}\partial
_{x}|z_{\varepsilon }|^{q}-q\mathcal{W}\partial _{x}\mathcal{W}%
|z_{\varepsilon }|^{q-1}\mathrm{sign}z_{\varepsilon }.
\end{align*}%
Integrating it in $x\in \mathbb{R}$, we have 
\begin{align*}
& \partial _{t}\int_{\mathbb{R}}|z_{\varepsilon }|^{q}\ dx\mathbf{\leqslant }%
C\left\{ \int_{\mathbb{R}}|\partial _{x}\mathrm{p}_{\varepsilon
}||z_{\varepsilon }|^{q-1}\ dx\right.  \\
& \left. +\Vert \partial _{x}\mathcal{W}\Vert _{L^{\infty }}\int_{\mathbb{R}%
}|z_{\varepsilon }|^{q}\ dx+\Vert \partial _{x}\mathcal{W}\Vert _{L^{\infty
}}\int_{\mathbb{R}}|\mathcal{W}||z_{\varepsilon }|^{q-1}\ dx\right\}  \\
& \mathbf{\leqslant }C\left\{ \Vert z_{\varepsilon }\Vert
_{L^{q}}^{q-1}+\Vert \partial _{x}\mathcal{W}\Vert _{L^{\infty }}\Vert
z_{\varepsilon }\Vert _{L^{q}}^{q}+\Vert \partial _{x}\mathcal{W}\Vert
_{L^{\infty }}\left( \Vert \mathcal{W}\Vert _{L^{q}}^{q}+\Vert
z_{\varepsilon }\Vert _{L^{q}}^{q}\right) \right\} ,
\end{align*}%
where we have used the estimate%
\begin{equation*}
\int_{\mathbb{R}}|\partial _{x}\mathrm{p}_{\varepsilon }||z_{\varepsilon
}|^{q-1}\ dx\mathbf{\leqslant }\Vert \partial _{x}\mathrm{p}_{\varepsilon
}\Vert _{L^{q}}\Vert z_{\varepsilon }\Vert _{L^{q}}^{q-1}\leqslant C\Vert
z_{\varepsilon }\Vert _{L^{q}}^{q-1}.
\end{equation*}%
following from \eqref{Z}$_{5}$. Therefore we obtain the Bihari type
inequality%
\begin{equation*}
\partial _{t}\Vert z_{\varepsilon }\Vert _{L^{q}}^{q}\mathbf{\leqslant }%
C\Vert z_{\varepsilon }\Vert _{L^{q}}^{q-1}+\alpha (t)\Vert z_{\varepsilon
}\Vert _{L^{q}}^{q}+\beta (t)
\end{equation*}%
with%
\begin{equation*}
\alpha (t)=C\Vert \partial _{x}\mathcal{W}\Vert _{L^{\infty }}\in
L^{1}(0,T),\qquad \beta (t)=C\Vert \partial _{x}\mathcal{W}\Vert _{L^{\infty
}}\Vert \mathcal{W}\Vert _{L^{q}}^{q}\in L^{1}(0,T),
\end{equation*}%
which gives%
\begin{equation*}
\Vert z_{\varepsilon }\Vert _{L^{q}}^{q}\leqslant C=C\left( \Vert
z_{\varepsilon }(0)\Vert _{L^{q}},\int_{0}^{T}\alpha (t)dt,\int_{0}^{T}\beta
(t)dt\right) \qquad \text{ for a.e. }t\in \lbrack 0,T].
\end{equation*}%
Using the assumptions (\ref{au}) we obtain the claimed estimate \eqref{Z}$%
_{3}$ for $z_{\varepsilon }$ in the $L^{q}$ norm.

The deduced estimates \eqref{Z}$_{3,4}$\ and the equation \eqref{Sys}$_{2}$
imply the last estimate \eqref{Z}$_{6},$ that ends the proof of this
lemma.\hfill $\blacksquare $

\bigskip

\bigskip

The estimates of Lemma \ref{lem1} do not imply a strong convergence of a
sub-sequence of $\left\{ z_{\varepsilon }\right\} ,$ that brings a
difficulty for the limit transition in the nonlinear terms of the system %
\eqref{Sys}. Hence in the sequel the principal aim is to derive the strong
convergence of $\left\{ z_{\varepsilon }\right\} .$ For it we use the
kinetic technique.

\bigskip \bigskip

\section{\protect\bigskip Limit transition by kinetic method}

\setcounter{equation}{0}

\bigskip

In this section, we apply the kinetic method in order to pass to the limit
the viscous system. \bigskip Let $\eta =\eta (s)\in C^{2}(\mathbb{R})$\ \ be
a convex function and consider the function $\ G,\ $such that $G^{\prime
}(s)=s\eta ^{\prime }(s).$ Multiplying the equation \eqref{Sys}$_{1}$ by $%
\eta ^{\prime }=\eta ^{\prime }(z_{\varepsilon })$ we show that for $P-$a.e. 
$\omega \in \Omega $ the pair $\eta (\cdot ,\omega )=\eta (z_{\varepsilon
})(\cdot ,\omega ),$ $\ G(\cdot ,\omega )=G(z_{\varepsilon }(\cdot ,\omega ))
$ satisfy the equality 
\begin{align}
\partial _{t}\eta +\partial _{x}G+\partial _{x}\mathrm{p}_{\varepsilon }\eta
^{\prime }& =\varepsilon \partial _{xx}\eta -\varepsilon (\partial
_{x}z_{\varepsilon })^{2}\eta ^{\prime \prime }-\partial _{x}\left( \eta 
\mathcal{W}\right)   \notag \\
& +\left( \eta -z_{\varepsilon }\eta ^{\prime }\right) \partial _{x}\mathcal{%
W}-\mathcal{W}\partial _{x}\mathcal{W}\eta ^{\prime }\qquad \text{in }%
\mathbb{R}_{T}\times \mathbb{R}.  \label{pe2}
\end{align}%
Hence for an arbitrary function $\varphi (t,x)\in C^{\infty
}([0,T];C_{0}^{\infty }(\mathbb{R}))$, such that $\varphi (T,\text{%
\textperiodcentered })=0,$\ and for $P-$a.e. $\omega \in \Omega $, the
stochastic process $z_{\varepsilon }(\cdot ,\omega )$ satisfies the entropy
type equality%
\begin{equation}
\begin{array}{l}
\displaystyle\int_{\mathbb{R}_{T}}\bigl\{\left( \eta \partial _{t}\varphi
+G\partial _{x}\varphi \right) -\left\{ \left( \partial _{x}\mathrm{p}%
_{\varepsilon }+\mathcal{W}\partial _{x}\mathcal{W}\right) \eta ^{\prime
}-\left( \eta -z_{\varepsilon }\eta ^{\prime }\right) \partial _{x}\mathcal{W%
}\right\} \varphi  \\ 
+\left( \eta \mathcal{W}-\varepsilon \,\partial _{x}\eta \right) \partial
_{x}\varphi \bigr\}dtdx+\displaystyle\int_{\mathbb{R}}\eta (z_{0,\varepsilon
})\varphi (0,\text{\textbf{\textperiodcentered }})\,dx=m_{\varepsilon ,\eta
}(\varphi )%
\end{array}
\label{2}
\end{equation}%
with%
\begin{equation*}
m_{\varepsilon ,\eta }(\varphi )=\int_{\mathbb{R}_{T}}\varepsilon (\partial
_{x}z_{\varepsilon })^{2}\eta ^{\prime \prime }(z_{\varepsilon })\varphi
\,dtdx
\end{equation*}%
such that $m_{\varepsilon ,\eta }(\varphi )\geqslant 0$ for any positive
function $\varphi .$

\bigskip Let us denote by%
\begin{equation*}
\left\vert s\right\vert _{\pm }=\max (\pm s,0),\qquad \mathrm{sign}_{+}(s)=%
\begin{cases}
1, & \text{if }s>0, \\ 
0, & \text{if }s\leqslant 0.%
\end{cases}%
\end{equation*}%
We introduce the convex function $\eta (z)=|z-c|_{+}$ and consider the
function 
\begin{equation*}
G(z)=\mathrm{sign}_{+}(z-c)\left( \frac{z^{2}}{2}-\frac{c^{2}}{2}\right) 
\text{\quad for any fixed parameter\quad }c\in \mathbb{R}.
\end{equation*}%
We notice that $\eta $ does not belong to $C^{2},$ then firstly we write (%
\ref{2}) for smooth convex functions $\eta _{\delta }(z)$\ and $G_{\delta
}(z)$, such that 
\begin{equation*}
\eta _{\delta }\rightarrow \eta \text{\quad strongly in\quad }%
W_{loc}^{1,\infty }(\mathbb{R})\text{ \quad as }\delta \rightarrow 0\text{%
\quad and\quad }G_{\delta }^{\prime }(s)=s\eta _{\delta }^{\prime }(s).
\end{equation*}%
Passing to the limit as $\delta \rightarrow 0$, we infer that (\ref{2})
holds for the above chosen pair $(\eta ,G).$ Since 
\begin{equation*}
\left( \eta -z_{\varepsilon }\eta ^{\prime }\right) \partial _{x}\mathcal{W}%
=-c\,\mathrm{sign}_{+}(z_{\varepsilon }-c)\partial _{x}\mathcal{W},
\end{equation*}%
we verify that for $P-$a.e. $\omega \in \Omega $, $z_{\varepsilon }(\cdot
,\omega )$\ satisfies the following equality 
\begin{equation}
\begin{array}{l}
\displaystyle\int_{\mathbb{R}_{T}}\biggl\{|z_{\varepsilon }-c|_{+}\partial
_{t}\varphi +\left( G(z_{\varepsilon })+|z_{\varepsilon }-c|_{+}\mathcal{W}%
\right) \partial _{x}\varphi \vspace{1mm} \\ 
-\left( \partial _{x}\mathrm{p}_{\varepsilon }+\mathcal{W}\partial _{x}%
\mathcal{W}+c\partial _{x}\mathcal{W}\right) \mathrm{sign}%
_{+}(z_{\varepsilon }-c)\varphi \vspace{1mm} \\ 
-\varepsilon \,\partial _{x}|z_{\varepsilon }-c|_{+}\partial _{x}\varphi %
\biggr\}dtdx\vspace{1mm} \\ 
+\displaystyle\int_{\mathbb{R}}|z_{0,\varepsilon }-c|_{+}\varphi (0,\text{%
\textbf{\textperiodcentered }})\,\,dx=m_{\varepsilon }(\varphi )%
\end{array}
\label{plus}
\end{equation}%
with 
\begin{equation}
m_{\varepsilon }(\varphi )=\int_{\mathbb{R}_{T}}\varepsilon \mathbf{\delta }%
(c=z_{\varepsilon }(t,x))|\nabla z_{\varepsilon }|^{2}\varphi \ dtdx.
\label{in}
\end{equation}%
Here $\mathbf{\delta }(s)$ denotes the Dirac function.

Therefore, for $P-$a.e. $\omega \in \Omega $, the stochastic process 
\begin{equation*}
f_{\varepsilon }(t,x,c)=\mathrm{sign}_{+}(z_{\varepsilon }(t,x)-c)
\end{equation*}%
\ solves the following integral problem%
\begin{equation}
\begin{array}{l}
\displaystyle\int_{\mathbb{R}_{T}}\biggl\{\int_{c}^{+\infty }f_{\varepsilon
}(t,x,s)\left( \partial _{t}\varphi +\left( s+\mathcal{W}\right) \partial
_{x}\varphi \right) \ ds\vspace{1mm} \\ 
-\left( \partial _{x}\mathrm{p}_{\varepsilon }+\mathcal{W}\partial _{x}%
\mathcal{W}+c\partial _{x}\mathcal{W}\right) f_{\varepsilon }(t,x,c)\varphi 
\vspace{1mm} \\ 
-\varepsilon f_{\varepsilon }(t,x,c)\partial _{x}z_{\varepsilon }\partial
_{x}\varphi \biggr\}\ dtdx\vspace{1mm} \\ 
+\int_{\mathbb{R}}|z_{0,\varepsilon }-c|_{+}\varphi (0,\cdot )\ dx\mathbf{=}%
m_{\varepsilon }(\varphi ).%
\end{array}
\label{in1}
\end{equation}%
Moreover defining the functional $L_{\varepsilon }$ as 
\begin{equation*}
L_{\varepsilon }(\psi )=\int_{\mathbb{R}}m_{\varepsilon }(\psi )\ dc\qquad 
\text{\ for any }\psi (t,x,c)\in C_{0}^{\infty }(\mathbb{R}_{T}\times 
\mathbb{R}),
\end{equation*}%
and using \eqref{in} and \eqref{in1}, we conclude that this functional is
positive and linear. Hence for $P-$a.e. $\omega \in \Omega $, by the Riesz
representation theorem (see for example, Theorem 1.39, p. 64 of \cite{EG}),
there exists an unique Radon measure $\mu _{\varepsilon }=\mu _{\varepsilon
}(t,x,c)$ defined on $\mathbb{R}_{T}\times \mathbb{R},$ such that 
\begin{equation*}
L_{\varepsilon }(\psi )=\int_{\mathbb{R}_{T}\times \mathbb{R}}\psi \ d\mu
_{\varepsilon }\qquad \text{\ for }\psi (t,x,c)\in C_{0}^{\infty }(\mathbb{R}%
_{T}\times \mathbb{R}).
\end{equation*}%
For $P-$a.e. $\omega \in \Omega $, the set of measures $\{\mu _{\varepsilon
}(\cdot ,\omega )\}$ is uniformly bounded on $\varepsilon ,$ since by %
\eqref{in}, \eqref{Z}$_{1}$ and for $P-$a.e. $\omega \in \Omega $, we have 
\begin{equation}
\int_{\mathbb{R}_{T}\times \mathbb{R}}\psi \ d\mu _{\varepsilon }=\int_{%
\mathbb{R}_{T}\times \mathbb{R}}\varepsilon \psi (t,x,z_{\varepsilon
}(t,x))|\nabla z_{\varepsilon }|^{2}\ dtdx<C  \label{eq2.55}
\end{equation}%
for any $\psi (t,x,c)\in C_{0}^{\infty }(\mathbb{R}_{T}\times \mathbb{R}),$
such that $0\leqslant \psi \leqslant 1.$

Analogously, if we take $\eta (z)=|z-c|_{-},\ \eta ^{\prime }(z)=\mathrm{sign%
}_{+}(z_{\varepsilon }-c)-1$\ \ and%
\begin{equation*}
G(z)=(\mathrm{sign}_{+}(z_{\varepsilon }-c)-1)\left( \frac{z^{2}}{2}-\frac{%
c^{2}}{2}\right) \text{\qquad for }c\in \mathbb{R}
\end{equation*}%
in (\ref{2}), and using a similar reasoning as in the deduction of (\ref{in1}%
), for $P-$a.e. $\omega \in \Omega $, we can show the following integral
relation 
\begin{equation}
\begin{array}{l}
\displaystyle\int_{\mathbb{R}_{T}}\biggl\{|z_{\varepsilon }-c|_{-}\,\partial
_{t}\varphi +\left( G(z_{\varepsilon })+|z_{\varepsilon }-c|_{-}\mathcal{W}%
\right) \partial _{x}\varphi \vspace{1mm} \\ 
-\left( \partial _{x}\mathrm{p}_{\varepsilon }+\mathcal{W}\partial _{x}%
\mathcal{W}+c\partial _{x}\mathcal{W}\right) (\mathrm{sign}%
_{+}(z_{\varepsilon }-c)-1)\varphi \vspace{1mm} \\ 
-\varepsilon \,\partial _{x}|z_{\varepsilon }-c|_{-}\partial _{x}\varphi %
\biggr\}{\ }dtdx\vspace{1mm} \\ 
+\displaystyle\int_{\mathbb{R}}|z_{0,\varepsilon }-c|_{-}\varphi (0,\text{%
\textbf{\textperiodcentered }})\,\,dx=m_{\varepsilon }(\varphi )%
\end{array}
\label{minus}
\end{equation}%
where $m_{\varepsilon }$ is defined by (\ref{in}). This last integral
relation can be written as%
\begin{equation}
\begin{array}{l}
\displaystyle\int_{\mathbb{R}_{T}}\biggl\{\int_{-\infty
}^{c}(1-f_{\varepsilon }(t,x,s))\{\partial _{t}\varphi +\left( s+\mathcal{W}%
\right) \partial _{x}\varphi \}\ ds\vspace{1mm} \\ 
-\left( \partial _{x}\mathrm{p}_{\varepsilon }+\mathcal{W}\partial _{x}%
\mathcal{W}+c\partial _{x}\mathcal{W}\right) (f_{\varepsilon
}(t,x,c)-1)\varphi \vspace{1mm} \\ 
-\varepsilon (f_{\varepsilon }(t,x,c)-1)\partial _{x}z_{\varepsilon
}\partial _{x}\varphi \biggr\}\ dtdx\vspace{1mm} \\ 
\displaystyle+\int_{\mathbb{R}}|z_{0,\varepsilon }-c|_{-}\varphi (0,\cdot
)dx=m_{\varepsilon }(\varphi ).%
\end{array}
\label{in2}
\end{equation}

Now, taking $\varphi =\partial _{c}\psi $, with $\psi (t,x,c)\in
C_{0}^{\infty }(\mathbb{R}_{T}\times \mathbb{R})$ in the equalities %
\eqref{in1} and (\ref{in2}), and next integrating by parts over $\mathbb{R}$
with respect to the parameter $c,$ we deduce that for $P-$a.e. $\omega \in
\Omega $, $f_{\varepsilon }$ satisfy the following two equations%
\begin{equation}
\begin{array}{l}
\displaystyle\int_{\mathbb{R}_{T}\times \mathbb{R}}\biggl\{f_{\varepsilon
}(t,x,c)\{\partial _{t}\psi +\left( c+\mathcal{W}\right) \partial _{x}\psi \}%
\vspace{1mm} \\ 
-\left( \partial _{x}\mathrm{p}_{\varepsilon }+\mathcal{W}\partial _{x}%
\mathcal{W}+c\partial _{x}\mathcal{W}\right) f_{\varepsilon }(t,x,c)\partial
_{c}\psi \vspace{1mm} \\ 
-\varepsilon f_{\varepsilon }(t,x,c)\partial _{x}z_{\varepsilon }\partial
_{x}\partial _{c}\psi \vspace{1mm}\biggr\}\ dtdxdc \\ 
\displaystyle+\int_{{\mathbb{R}}^{2}}\mathrm{sign}_{+}\left(
z_{0,\varepsilon }-c\right) \psi (0,\cdot )\ dxdc\mathbf{=}\int_{\mathbb{R}%
_{T}\times \mathbb{R}}\partial _{c}\psi \ d\mu _{\varepsilon }%
\end{array}
\label{f0}
\end{equation}%
and 
\begin{equation}
\begin{array}{l}
\displaystyle\int_{\mathbb{R}_{T}\times \mathbb{R}}\biggl\{(f_{\varepsilon
}(t,x,c)-1)\{\partial _{t}\psi +\left( c+\mathcal{W}\right) \partial
_{x}\psi \}\vspace{1mm} \\ 
-\left( \partial _{x}\mathrm{p}_{\varepsilon }+\mathcal{W}\partial c\mathcal{%
W}+c\partial _{x}\mathcal{W}\right) (f_{\varepsilon }(t,x,c)-1)\partial
_{c}\psi \vspace{1mm} \\ 
-\varepsilon (f_{\varepsilon }(t,x,c)-1)\ \partial _{x}\partial _{c}\psi %
\biggr\}\ dtdxdc \\ 
\displaystyle+\int_{\mathbb{R}^{2}}\left( \mathrm{sign}_{+}\left(
z_{0,\varepsilon }-c\right) -1\right) \psi (0,\cdot )\ dx=\int_{\mathbb{R}%
_{T}\times \mathbb{R}}\partial _{c}\psi \ d\mu _{\varepsilon }.%
\end{array}
\label{f1}
\end{equation}

\bigskip We emphasize that due to the presence of the term $\partial _{x}%
\mathrm{p}_{\varepsilon }$, the equations \eqref{f1}-\eqref{f1} are
nonlinear, requiring a strong convergence of the sequence $\{f_{\varepsilon
}\}$ in order to pass to the limit, as $\varepsilon \rightarrow 0$. The next
lemma establish helpful estimates for this analysis.

\begin{lemma}
\label{Lemma 4.2} For $P-$a.e. $\omega\in\Omega$\ there exist some constants 
$C(n),$ which are independent of $\varepsilon $, but may depend of $\omega $
and of an arbitrary fixed natural $n>0$, such that 
\begin{eqnarray}
0 &\leqslant &f_{\varepsilon }\leqslant 1\quad \text{ a.e. in }\mathbb{R}%
_{T}\times \mathbb{R},\quad \quad \text{ }\frac{\partial f_{\varepsilon }}{%
\partial c}\leqslant 0\quad \text{ in \ }\mathcal{D}^{\prime }(\mathbb{R}%
_{T}\times \mathbb{R}),  \notag \\
z_{\varepsilon }(t,x) &=&\int_{0}^{+\infty }f_{\varepsilon
}(t,x,s)\,ds-\int_{-\infty }^{0}(1-f_{\varepsilon }(t,x,s))\,ds  \label{ff0}
\end{eqnarray}%
{and 
\begin{equation}
||\partial _{t}f_{\varepsilon }||_{L^{2}(0,T;\,H^{-2}(Y))}\leqslant C(n),
\label{timedelta}
\end{equation}%
} where $Y=(-n,n)\times (-n,n)$.
\end{lemma}

\textbf{Proof}. The relations \eqref{ff0} follow directly from the
properties of the stochastic processes $f_{\varepsilon }(t,x,c)=\mathrm{sign}%
_{+}(z_{\varepsilon }(t,x)-c),$ which decreases on $c\in \mathbb{R}.$

Applying Holder's inequality and the continuous embedding $%
H_{0}^{2}(Y)\hookrightarrow C_{0}(Y)$, the relation \eqref{f0} yields%
\begin{eqnarray*}
|(\partial _{t}f_{\varepsilon }\mathbf{,}\psi )| &\leqslant &C\left( \Vert c+%
\mathcal{W}\Vert _{L^{2}(Y)}+||\partial _{x}\mathrm{p}_{\varepsilon }+%
\mathcal{W}\partial _{x}\mathcal{W}+c\partial _{x}\mathcal{W}%
||_{L^{2}(Y)}\right) \left\Vert \nabla _{x,c}\psi \vspace{1mm}\right\Vert
_{L^{2}(Y)} \\
&&+\varepsilon \left\Vert \partial _{x}z_{\varepsilon }\right\Vert
_{2}\left\Vert \partial _{x}\partial _{c}\psi \vspace{1mm}\right\Vert
_{L^{2}(Y)}+\int_{(0,T)\times Y}1\ d\mu _{\varepsilon }\Vert \partial
_{c}\psi \Vert _{L^{\infty }(Y)} \\
&\leqslant &C(n){g(t)}\Vert \psi \Vert _{H^{2}({\mathbb{R}}^{2})}\qquad 
\text{\ for any }\psi \in H_{0}^{2}(Y),
\end{eqnarray*}%
with ${g(t)=}1+\Vert \mathcal{W}\Vert _{L^{2}(Y)}+||\mathcal{W}\partial _{x}%
\mathcal{W}||_{L^{2}(Y)}+||\partial _{x}\mathcal{W}||_{L^{2}(Y)}+\left\Vert 
\sqrt{\varepsilon }\partial _{x}z_{\varepsilon }\right\Vert _{L^{2}(Y)}\in
L^{1}(0,T)$ by \eqref{Z}$.$ Therefore%
\begin{equation*}
||\partial _{t}f_{\varepsilon }||_{H^{-2}(Y)}=\sup_{\psi \in
H_{0}^{2}(Y)}\left\{ |(\partial _{t}f_{\varepsilon },\psi )|:\quad ||\psi
||_{H^{2}(Y)}\leqslant 1\right\} \leqslant C(n){g}_{{\varepsilon }}{(t)}
\end{equation*}%
Due to \eqref{au} and \eqref{Z}$_{1}$ we have $||{g}_{{\varepsilon }%
}||_{L^{2}(0,T)}\leqslant C,$ that gives the estimate (\ref{timedelta}).$%
\hfill \;\blacksquare $

\bigskip

\bigskip

\subsection{Passage to the limit on the nonlinear kinetic equations}

\label{sec32}

\bigskip

The estimates \eqref{limit0}$_{2}$ \eqref{Z}, \eqref{eq2.55} and \eqref{ff0}
imply that for $P-$a.e. $\omega \in \Omega $ there exist sub-sequences of $%
\{z_{\varepsilon },f_{\varepsilon },\mathrm{p}_{\varepsilon },\mu
_{\varepsilon }\},$ which may depend on $\omega $, such that 
\begin{align}
\varepsilon \nabla z_{\varepsilon }& \rightarrow 0\quad \text{strongly in $%
L^{2}(\mathbb{R}_{T})$},  \notag \\
z_{\varepsilon }& \rightharpoonup z\quad \text{weakly -- $\ast $ in $%
L^{\infty }(0,T;L^{2}(\mathbb{R})\cap L^{q}(\mathbb{R})),$}  \notag \\
f_{\varepsilon }& \rightharpoonup f\quad \text{weakly -- $\ast $ in }%
L^{\infty }(\mathbb{R}_{T}\times \mathbb{R}\text{$),$}  \notag \\
\mathrm{p}_{\varepsilon },\ \partial _{x}\mathrm{p}_{\varepsilon }&
\rightharpoonup f,\ \mathrm{p},\ \partial _{x}\mathrm{p}\quad \text{weakly
-- $\ast $ in }L^{\infty }(0,T;L^{p}(\mathbb{R}),\qquad \forall p\in \left[
1,+\infty \right] ,  \notag \\
\partial _{xx}^{2}\mathrm{p}_{\varepsilon }& \rightharpoonup \partial
_{xx}^{2}\mathrm{p}_{\varepsilon }\quad \text{weakly -- $\ast $ in }%
L^{\infty }(0,T;L^{q/2}(\mathbb{R}))\text{$,$}  \notag \\
\mu _{\varepsilon }& \rightharpoonup \mu \quad \text{weakly in }\mathcal{M}%
_{loc}^{+}(\mathbb{R}_{T}\times \mathbb{R}).  \label{limit0}
\end{align}%
In the deduction of (\ref{limit0})$_{6}$, we applied the Theorem 1.41, p. 66
of \cite{EG}.

In order to deal with the limit of the nonlinear term $\left( \partial _{x}%
\mathrm{p}_{\varepsilon }\right) f_{\varepsilon }$\ in the equations %
\eqref{f0}-\eqref{f1}, we first write a lemma which consists in an
adaptation to our setting of the Theorem 6, p. 7 of \cite{Ev}.

\begin{lemma}
\label{16_09} Let $Y$ be a bounded open subset of $\mathbb{R}^2.$ For any $%
r>1$, the following compact embedding holds 
\begin{equation*}
L^{\infty }(Y)\hookrightarrow (W^{1,r}(Y))^*.
\end{equation*}
\end{lemma}

\textbf{Proof.} Let $\{g_{n}\}$ be a bounded sequence in $L^{\infty }(Y)$,
and denote by $M=\sup_{n}\Vert g_{n}\Vert _{L^{\infty }(Y)}.$ We can
consider $\{g_{n}\}$ as a bounded sequence of functionals in $%
(W^{1,r}(Y))^{\ast },$ such that 
\begin{equation*}
|\langle g_{n},\phi \rangle |=\left\vert \int_{Y}\phi g_{n}dx\right\vert
\leqslant MC(Y,r)\Vert \phi \Vert _{L^{r}(Y)}\leqslant MC(Y,r)\Vert \phi
\Vert _{W^{1,r}(Y)},\quad \forall \phi \in W^{1,r}(Y).
\end{equation*}%
Since $||g_{n}||_{L^{\infty }(Y)}\leqslant M,$ then there exists a
sub-sequence $\{g_{n_{k}}\}$ such that 
\begin{equation*}
g_{n_{k}}\rightharpoonup g\quad \text{weakly in }\quad L^{\infty }(Y),
\end{equation*}%
that implies 
\begin{equation}
g_{n_{k}}\rightharpoonup g\quad \text{weakly in }\quad (W^{1,r}(Y))^{\ast }.
\label{16_6}
\end{equation}%
Let us denote by $B$ the closed unit ball in $W^{1,r}(Y)$. From Theorems 5,
p. 7 of \cite{Ev}, $B$ is a compact set in $L^{p}(Y)$ with $1\leqslant p<%
\frac{2r}{2-r}$. Given $\delta >0$, we can take $\phi _{i}\in B$, $i=1,\dots
N(\delta )$, such that for any $\phi \in B$, there exists some $i\in
\{1,\dots ,N(\delta )\}$ that verifies 
\begin{equation}
\Vert \phi -\phi _{i}\Vert _{L^{p}(Y)}\leqslant \frac{\delta }{4MC(Y,p)}.
\label{16_5}
\end{equation}

On the other hand, due to the weak convergence \eqref{16_6}, there exists $%
N\in \mathbb{N}$ such that $\forall k>N$, we have 
\begin{equation*}
\left\vert \int_{Y}\phi _{i}g_{n_{k}}dx-\int_{Y}\phi _{i}gdx\right\vert
\leqslant \frac{\delta }{2},\quad \forall i=1,2\dots ,N(\delta ).
\end{equation*}%
For any $\phi \in B,$ choosing $\phi _{i}$ verifying \eqref{16_5}, we deduce 
\begin{equation*}
|\langle g_{n_{k}}-g,\phi \rangle |\leqslant 2MC(Y,p)\Vert \phi -\phi
_{i}\Vert _{L^{p}(Y)}+\left\vert \int_{Y}\phi _{i}g_{n_{k}}dx-\int_{Y}\phi
_{i}gdx\right\vert \leqslant \delta .
\end{equation*}%
Therefore, $\forall k>N$ we obtain 
\begin{equation*}
\Vert g_{n_{k}}-g\Vert _{(W^{1,r}(Y))^{\ast }}=\sup_{\phi \in B}|\langle
g_{n_{k}}-g,\phi \rangle |\leqslant \delta ,
\end{equation*}%
which means that $g_{n_{k}}\rightarrow g$ strongly in $(W^{1,r}(Y))^{\ast }.$

$\hfill \;\blacksquare$

\bigskip

For the $q$ introduced in \eqref{au}, we have 
\begin{equation*}
L^{\infty }(Y)\subseteq (W^{1,q/2}(Y))^{\ast }\subseteq H^{-2}(Y),
\end{equation*}%
where, accordingly to Lemma \ref{16_09}, $L^{\infty }(Y)$ is compactly
embedded in $(W^{1,q/2}(Y))^{\ast }$, and in addition $(W^{1,q/2}(Y))^{\ast
} $ is continuously embedded in $H^{-2}(Y)$.

For $P-$a.e. $\omega \in \Omega $, the sequence $\{f_{\varepsilon }(t,x,c)\}$
is bounded in $L^{\infty }(\mathbb{R}_{T}\times \mathbb{R})$, it is also
bounded in $L^{\infty }(0,T;L^{\infty }(Y))$ with $Y=(-n,n)\times (-n,n)$,
for any $n\in \mathbb{N}$. Hence {the estimate (\ref{timedelta})} {and} the
Aubin-Lions-Simon compactness result \cite{sim}, (see also, Theorem II.5.16,
p. 102 of \cite{Boyer}) imply that there exists a sub-sequence of \ $%
\{f_{\varepsilon }\}$ (after relabeling by the same index $\varepsilon )$
such that%
\begin{equation*}
f_{\varepsilon }\rightarrow f\quad \text{strongly in $L^{\infty
}(0,T;(W^{1,q/2}(Y))^{\ast })$}
\end{equation*}%
for each fixed $n>0.$ Therefore a diagonalization procedure shows that there
exists a subsequence of $\{f_{\varepsilon }\},$ such that%
\begin{equation*}
f_{\varepsilon }\rightarrow f\quad \text{strongly in $L^{\infty
}(0,T;(W_{loc}^{1,q/2}(\mathbb{R}^{2}))^{\ast }).$}
\end{equation*}%
The a priori estimate \eqref{Z}$_{7,8}$ give 
\begin{equation}
\partial _{x}\mathrm{p}_{\varepsilon }\rightharpoonup \partial _{x}\mathrm{p}%
\quad \text{weakly-}\ast \text{ in }L^{\infty }(0,T;W_{loc}^{1,q/2}(\mathbb{R%
}^{2})).  \label{op}
\end{equation}%
Combining these two convergences, we obtain 
\begin{equation*}
\int_{\mathbb{R}_{T}\times \mathbb{R}}\left( \partial _{x}\mathrm{p}%
_{\varepsilon }\right) f_{\varepsilon }\psi \ dtdxdc\rightarrow \,\int_{%
\mathbb{R}_{T}\times \mathbb{R}}\left( \partial _{x}\mathrm{p}\right) f\psi
\ dtdxdc
\end{equation*}%
for any smooth function $\psi $ having compact support on $\mathbb{R}%
_{T}\times \mathbb{R}.$ Therefore the equalities \eqref{f0} and \eqref{f1}
implies that $z,f,\mathrm{p}$ and $\mu ,$ obtained by the limit transition %
\eqref{limit0}, satisfy the integral relations%
\begin{eqnarray}
&&\int_{\mathbb{R}_{T}\times \mathbb{R}}\bigl\{f(t,x,c)\{\partial _{t}\psi
+\left( c+\mathcal{W}\right) \partial _{x}\psi \}\vspace{1mm}  \notag \\
&&-\left( \partial _{x}\mathrm{p}+\mathcal{W}\partial _{x}\mathcal{W}%
+c\partial _{x}\mathcal{W}\right) f(t,x,c)\partial _{c}\psi \bigr\}\ dtdxdc 
\notag \\
&&+\int_{{\mathbb{R}}^{2}}\mathrm{sign}_{+}\left( z_{0}-c\right) \psi
(0,\cdot )\ dxdc\mathbf{=}\int_{\mathbb{R}_{T}\times \mathbb{R}}\partial
_{c}\psi \,\,\ d\mu .  \label{2.6}
\end{eqnarray}%
and%
\begin{eqnarray}
&&\int_{\mathbb{R}_{T}\times \mathbb{R}}\bigl\{(f(t,x,c)-1)\{\partial
_{t}\psi +\left( c+\mathcal{W}\right) \partial _{x}\psi \}\vspace{1mm} 
\notag \\
&&-\left( \partial _{x}\mathrm{p}+\mathcal{W}\partial _{x}\mathcal{W}%
+c\partial _{x}\mathcal{W}\right) (f(t,x,c)-1)\partial _{c}\psi \bigr\}\
dtdxdc  \notag \\
&&+\int_{\mathbb{R}^{2}}\left( \mathrm{sign}_{+}\left( z_{0}-c\right)
-1\right) \psi (0,\cdot )\ dx=\int_{\mathbb{R}_{T}\times \mathbb{R}}\partial
_{c}\psi \,\,\ d\mu ,  \label{2.7}
\end{eqnarray}%
where $\psi $ is a smooth function having compact support on $\mathbb{R}%
_{T}\times \mathbb{R}.$

\bigskip

\bigskip

\begin{lemma}
\label{Lemma 4.44} For a.s. $\omega \in \Omega $, we have 
\begin{equation}
0\leqslant f\leqslant 1\quad \text{a.e. in }\mathbb{R}_{T}\times \mathbb{R}%
\quad \text{ and}\quad \text{ }\frac{\partial f}{\partial c}\leqslant
0\qquad \text{in \ }\mathcal{D}^{\prime }(\mathbb{R}_{T}\times \mathbb{R}).
\label{fff}
\end{equation}%
Moreover, for a.s. $\omega \in \Omega $, the limit stochastic processes 
\begin{equation*}
z\in L^{\infty }(0,T;L^{2}(\mathbb{R})\cap L^{q}(\mathbb{R}))\quad \text{and}%
\quad \widehat{z}\in L^{\infty }(0,T;L^{q/p}(\mathbb{R}))\qquad \forall p\in
\lbrack 1,q),
\end{equation*}%
defined in $(\ref{limit0})_{2}$ and \eqref{vvv}, respectively, verify 
\begin{eqnarray}
z(t,x) &=&\int_{0}^{+\infty }f(t,x,c)\,dc-\int_{-\infty
}^{0}(1-f(t,x,c))\,dc,  \notag \\
\widehat{z}(t,x) &=&\int_{0}^{+\infty }f(t,x,c)p{c}^{p-1}dc+\int_{-\infty
}^{0}(1-f(t,x,c))p|c|^{p-1}dc.  \label{eq2.88}
\end{eqnarray}
\end{lemma}

\textbf{Proof}. The relations \eqref{fff} follow from \eqref{ff0}$_{1}$ and
the weak convergences \eqref{limit0}$_{3}$.

Now let us consider an arbitrary $p\in \lbrack 1,q)$ and a positive test
function $\varphi =\varphi (t,x)\in \mathcal{D}(\mathbb{R}_{T}).$ We have
the estimate 
\begin{eqnarray*}
\int_{\mathbb{R}_{T}}\int_{0}^{+\infty }\varphi f_{\varepsilon
}pc^{p-1}dtdxdc &=&\int_{\mathbb{R}_{T}}\varphi |z_{\varepsilon
}|_{+}^{p}dtdx \\
&\leq &||z_{\varepsilon }||_{L^{\infty }(0,T;L^{q}(\mathbb{R}))}^{p}\Vert
\varphi \Vert _{L^{1}(0,T;L^{r}(\mathbb{R}))}\leq C(p,\varphi )
\end{eqnarray*}%
with $r,$ such that $1=\frac{p}{q}+\frac{1}{r},$ and the constant $%
C(p,\varphi )$ depending only on $p$ and $\varphi $.

By the weak convergence of $f_{\varepsilon }$ to $f$ in $L^{\infty }(\mathbb{%
R}_{T}\times \mathbb{R}),$ \ for any $M>0$ and any $\delta >0$, there exists 
$\varepsilon _{1}=\varepsilon _{1}(M,\delta )>0,$ such that for each
positive $\varepsilon <\varepsilon _{1}$, we have 
\begin{equation}
\left\vert \int_{\mathbb{R}_{T}}\int_{0}^{M}\varphi f_{\varepsilon
}pc^{p-1}dtdxdc-\int_{\mathbb{R}_{T}}\int_{0}^{M}\varphi
fpc^{p-1}dtdxdc\right\vert \leq \frac{\delta }{3}.  \label{aa66}
\end{equation}%
If we take $\delta =1$, we obtain that 
\begin{equation*}
\int_{\mathbb{R}_{T}}\int_{0}^{M}\varphi fpc^{p-1}dtdxdc\leq \frac{1}{3}%
+\int_{\mathbb{R}_{T}}\int_{0}^{+\infty }\varphi f_{\varepsilon
}pc^{p-1}dtdxdc\leq \frac{1}{3}+C(p,\varphi )\quad \text{for all \ }M,
\end{equation*}%
which implies 
\begin{equation}
\int_{\mathbb{R}_{T}}\int_{0}^{+\infty }\varphi fpc^{p-1}dtdtdc<\infty
\label{aa44}
\end{equation}%
by the monotone convergence theorem.

For any $M>0$, using the H\"{o}lder inequality and the estimate (\ref{Z})$_{%
\text{2,3}},$ we verify that 
\begin{equation}
\left\vert \int_{\mathbb{R}_{T}}\int_{M}^{+\infty }\varphi f_{\varepsilon
}pc^{p-1}dtdxdc\right\vert \leq \frac{C\Vert \varphi \Vert _{L^{\infty }(%
\mathbb{R}_{T})}}{M^{q-p}}\int_{|z_{\varepsilon }|_{+}\geq M}|z_{\varepsilon
}|_{+}^{q}dtdx\leq \frac{C(\varphi ,p)}{M^{q-p}}.  \label{aa11}
\end{equation}%
With the help of the inequalities \eqref{aa66}-\eqref{aa11} and considering $%
M$ large enough, we deduce 
\begin{equation*}
\left\vert \int_{\mathbb{R}_{T}}\int_{0}^{+\infty }\varphi f_{\varepsilon
}pc^{p-1}dtdxdc-\int_{\mathbb{R}_{T}}\int_{0}^{+\infty }\varphi
fpc^{p-1}dtdtdc\right\vert \leq \delta ,
\end{equation*}%
which corresponds to 
\begin{equation}
|z_{\varepsilon }|_{+}^{p}=\int_{0}^{+\infty }f_{\varepsilon }(\mathbf{\cdot 
},\mathbf{\cdot },c)pc^{p-1}\ dc\rightharpoonup \int_{0}^{+\infty }f(\mathbf{%
\cdot },\mathbf{\cdot },c)pc^{p-1}\ dc\quad \text{in \ }\mathcal{D}^{\prime
}(\mathbb{R}_{T}),  \label{a0}
\end{equation}%
Dealing with the function $1-f_{\varepsilon }$ instead of $f_{\varepsilon }$%
, and applying the same reasoning, we can derive 
\begin{equation}
|z_{\varepsilon }|_{-}^{p}=\int_{-\infty }^{0}(1-f_{\varepsilon }(\mathbf{%
\cdot },\mathbf{\cdot },c))p|c|^{p-1}\ dc\rightharpoonup \int_{-\infty
}^{0}(1-f(\mathbf{\cdot },\mathbf{\cdot },c))p|c|^{p-1}\ dc\quad \text{in \ }%
\mathcal{D}^{\prime }(\mathbb{R}_{T}).  \label{a00}
\end{equation}%
\hfill

Now let us consider 2 distinct cases:

1. Let us consider $p=1.$ Having the weak convergence of $z_{\varepsilon
}=|z_{\varepsilon }|_{+}-|z_{\varepsilon }|_{+}$ to $z$ in $L^{\infty
}(0,T;L^{2}(\mathbb{R})\cap L^{q}(\mathbb{R}))$ by \eqref{limit0}$_{2},$ we
see that the convergences \eqref{a0}-\eqref{a00} imply \eqref{eq2.88}$_{1}$.

2. Now let us consider that $p\in \lbrack 1,q)$. Due to estimates (\ref{Z})$%
_{\text{3}}$ the set 
\begin{equation}
\{|z_{\varepsilon }|^{p},\ |z_{\varepsilon }|_{+}^{p},\ |z_{\varepsilon
}|_{-}^{p}\}_{\varepsilon >0}\text{ is uniformly bounded in }L^{\infty
}(0,T;L^{q/p}(\mathbb{R})),  \label{es}
\end{equation}%
therefore we have 
\begin{equation}
|z_{\varepsilon }|^{p}=|z_{\varepsilon }|_{+}^{p}+|z_{\varepsilon
}|_{-}^{p}\rightharpoonup \widehat{z}\quad \text{weakly -- $\ast $ in $%
L^{\infty }(0,T;L^{q/p}(\mathbb{R})).$}  \label{vvv}
\end{equation}%
The limit function $\widehat{z}\in L^{\infty }(0,T;L^{q/p}(\mathbb{R}))$ has
the representation \eqref{eq2.88}$_{2}$ by \eqref{a0}-\eqref{a00}.\hfill $%
\blacksquare $

\bigskip

In the sequel we will need the following properties of behavior at infinity.

\begin{lemma}
\label{lemma22} For $P-$a.e. $\omega\in\Omega$, a.e. $t\in \lbrack 0,T]$ and
any $p\in \lbrack 1,q)$ we have 
\begin{eqnarray}
&&\lim_{c\rightarrow +\infty }c^{p}\int_{\mathbb{R}}f(t,x,c)dx=0,\qquad
\lim_{c\rightarrow -\infty }|c|^{p}\int_{\mathbb{R}}(1-f(t,x,c))d{x}=0,
\label{equ1} \\
&&\lim_{c\rightarrow \pm \infty }\mu (c)=0.  \label{equ2}
\end{eqnarray}
\end{lemma}

\noindent \textbf{Proof.} By \eqref{eq2.88}$_{2}$ we have that for a.e. $%
t\in \lbrack 0,T]$%
\begin{eqnarray*}
\infty &>&\int_{\mathbb{R}}\left[ \int_{0}^{+\infty }f(t,x,c){c}%
^{p-1}dc+\int_{-\infty }^{0}(1-f(t,x,c))|c|^{p-1}ds\right] d{x} \\
&&\left[ \int_{0}^{+\infty }\left( \int_{\mathbb{R}}f(t,x,c)dx\right) {c}%
^{p-1}dc+\int_{-\infty }^{0}\left( \int_{\mathbb{R}}(1-f(t,x,c))d{x}\right)
|c|^{p-1}ds\right] ,
\end{eqnarray*}%
which implies (\ref{equ1}), knowing that $0\leqslant f\leqslant 1$.

Finally, to show (\ref{equ2}) we use equalities (\ref{2.6}) and (\ref{2.7}),
with the help of (\ref{equ1}).\hfill $\blacksquare $

\bigskip \bigskip

\subsection{Strong convergence}

\bigskip

The integral equalities (\ref{2.6}) and (\ref{2.7}) can be written as the
following transport type equation for $g=f$ and $g=f-1$ \ in the
distribution sense in $\mathcal{D}^{\prime }(\mathbb{R}_{T}\times \mathbb{R})
$ 
\begin{equation}
g_{t}+\partial _{x}\left( g\left( c+\mathcal{W}\right) \right) -\partial
_{c}\left( g\left( \partial _{x}\mathrm{p}+\mathcal{W}\partial _{x}\mathcal{W%
}+c\partial _{x}\mathcal{W}\right) \right) \mathbf{=}\partial _{c}\mu ,
\label{eq2.9}
\end{equation}%
that yields the next result.

\begin{lemma}
\label{lemma4} For $P-$a.e. $\omega \in \Omega $, the function $F=f(1-f)$
satisfies%
\begin{equation}
\partial _{t}F+\partial _{x}\left( F\left( c+\mathcal{W}\right) \right)
-\partial _{c}\left( F\left( \partial _{x}\mathrm{p}+\mathcal{W}\partial _{x}%
\mathcal{W}+c\partial _{x}\mathcal{W}\right) \right) \leqslant 0\quad \text{
in \ }\mathcal{D}^{\prime }(\mathbb{R}_{T}\times \mathbb{R}).  \label{RRR}
\end{equation}
\end{lemma}

\noindent \textbf{Proof. }Let $\rho _{\delta }=\rho _{\delta }(s)$ be the
standard symmetric kernel on the variable $s\mathbf{\in }\mathbb{R}$. We
denote by 
\begin{equation*}
u^{\delta }(x,c)=\int_{\mathbb{R}^{2}}u(s,y)\rho _{\delta }(x-y)\rho
_{\delta }(c-s)\ dyds
\end{equation*}%
the mollification of a function $u=u(x,c)$ on both variables $c,$ $x$.

Regularizing the equation (\ref{eq2.9}) for $g=f$ by this mollification
operator, we obtain the following equality%
\begin{equation}
\partial _{t}f^{\delta }+\partial _{x}\left( f^{\delta }\left( c+\mathcal{W}%
\right) \right) -\partial _{c}\left( f^{\delta }\left( \partial _{x}\mathrm{p%
}+\mathcal{W}\partial _{x}\mathcal{W}+c\partial _{x}\mathcal{W}\right)
\right) \mathbf{=}\partial _{c}[\,\mu ]^{\delta }+\epsilon _{1}^{\delta }
\label{R1}
\end{equation}%
in\ $\mathcal{D}^{\prime }(\mathbb{R}_{T}\times \mathbb{R}),$ with the
remainder%
\begin{eqnarray*}
\epsilon _{1}^{\delta } &=&\partial _{x}\left( f^{\delta }\left( c+\mathcal{W%
}\right) -\left[ f\left( c+\mathcal{W}\right) \right] ^{\delta }\right)  \\
&&-\partial _{c}\left( f^{\delta }\left( \partial _{x}\mathrm{p}+c\partial
_{x}\mathcal{W}+\mathcal{W}\partial _{x}\mathcal{W}\right) -\left[ f\left(
\partial _{x}\mathrm{p}+c\partial _{x}\mathcal{W}+\mathcal{W}\partial _{x}%
\mathcal{W}\right) \right] ^{\delta }\right) .
\end{eqnarray*}%
Since $f\in L^{\infty }(\mathbb{R}_{T}\times \mathbb{R}),\ \ \partial _{x}%
\mathrm{p}\in L^{\infty }(0,T;W_{loc}^{1,q/2}(\mathbb{R})),$ $c\in
W_{loc}^{1,\infty }(\mathbb{R})\ $and $\mathcal{W}$ fulfills (\ref{au})\ by
Lemma 2.2, p. 972-06, \ \cite{lellis}, we derive that 
\begin{equation*}
\lim_{\delta \rightarrow 0}\epsilon _{1}^{\delta }=0\text{\qquad strongly
in\qquad }L_{loc}^{1}(\mathbb{R}_{T}\times \mathbb{R}).
\end{equation*}

By the similar considerations we deduce that the equation (\ref{eq2.9}) for $%
g=f-1$ gives%
\begin{eqnarray}
&&(f^{\delta }-1)_{t}+\partial _{x}\left( (f^{\delta }-1)\left( c+\mathcal{W}%
\right) \right) -\partial _{c}\left( (f^{\delta }-1)\left( \partial _{x}%
\mathrm{p}+c\partial _{x}\mathcal{W}+\mathcal{W}\partial _{x}\mathcal{W}%
\right) \right)   \notag \\
&=&\partial _{c}\,[\,\mu ]^{\delta }+\epsilon _{2}^{\delta }\qquad \text{in}%
\;\;\mathcal{D}^{\prime }(\mathbb{R}_{T}\times \mathbb{R}),  \label{R2}
\end{eqnarray}%
with%
\begin{equation*}
\lim_{\delta \rightarrow 0}\epsilon _{2}^{\delta }=0\text{\qquad strongly
in\qquad }L_{loc}^{1}(\mathbb{R}_{T}\times \mathbb{R}).
\end{equation*}

Since $f^{\delta }$ is smooth in $(x,c)$, by multiplying the equality %
\eqref{R1} by $(f^{\delta }-1)$ and the equality \eqref{R2} by $f^{\delta }$
and taking the limit transitions on $\delta \rightarrow 0,$ we deduce that $%
F=-f(f-1)$ satisfies (\ref{RRR}), where we have used the following relation 
\begin{equation*}
\int_{\mathbb{R}}\partial _{c}[\,\mu ]^{\delta ,\theta }\ (f^{\delta ,\theta
}-1)+\partial _{c}[\,\mu ]^{\delta ,\theta }\ f^{\delta ,\theta }\
\,dc=-\int_{\mathbb{R}}\partial _{c}f^{\theta ,\delta }\,2\left[ \mu \right]
^{\delta ,\theta }dc\geqslant 0.
\end{equation*}%
The latter inequality follows from properties \eqref{eq2.88} and $\mu $ is a
positive measure.

\hfill $\blacksquare $

\bigskip

The equalities (\ref{2.6}) and (\ref{2.7}) yield, for $P-$a.e. $%
\omega\in\Omega$, the existence of traces for the functions $f$ and $1-f$ at
the time $t=0$.

\begin{lemma}
\label{lemma2} For $P-$a.e. $\omega\in\Omega$, the function $f$ has the
trace $f^{0}=f^{0}(x,c)$ at the time moment $t=0,$ such that for $P-$a.e. $%
\omega\in\Omega$ 
\begin{equation}
f^{0}=\lim_{\delta \rightarrow 0^{+}}\frac{1}{\delta }\int_{0}^{\delta
}f(t,\cdot ,\mathbf{\cdot })\,dt\mathbb{\quad }\text{and}\mathbb{\quad }%
f^{0}-\left( f^{0}\right) ^{2}=0\text{\ \ a.e. on }{\mathbb{R}}^{2}.
\label{2.13}
\end{equation}
\end{lemma}

\noindent \textbf{Proof.} For $P-$a.e. $\omega \in \Omega $, let us
introduce the vector function%
\begin{equation*}
\mathbf{F}(t,x,c)=\left( f,\ f\left( c+\mathcal{W}\right) ,\ -\left( f\left(
\partial _{x}\mathrm{p}+\mathcal{W}\partial _{x}\mathcal{W}+c\partial _{x}%
\mathcal{W}\right) +\mu \right) \right) 
\end{equation*}%
which is divergence free in $\mathcal{D}^{\prime }(\mathbb{R}_{T}\times 
\mathbb{R})$ by the equation \eqref{eq2.9} and belongs to the space $%
\mathcal{\mathcal{DM}}^{ext}(\mathbb{R}_{T}\times \mathbb{R})$ of extended
divergence-measure field over $\mathbb{R}_{T}\times \mathbb{R}$ in
accordance of Definition 1.1 of \cite{chen} by (\ref{au}), (\ref{Z}).
Therefore\ it\ follows from Theorem 2.2 in \cite{chen} that $\mathbf{F}\cdot
(1,0,0)$\ at the time moment $t=0$ can be considered as a continuous linear
functional defined over the space $C_{0}^{1}({\mathbb{R}}^{2}),$ such that 
\begin{equation*}
\left\langle \mathbf{F}\cdot (1,0,0)|_{t=0},\phi \right\rangle =\lim_{\delta
\rightarrow 0^{+}}\int_{{\mathbb{R}}^{2}}\left( \frac{1}{\delta }%
\int_{0}^{\delta }f(t,x,c)\,dt\right) \phi (x,c)\,\,dxdt
\end{equation*}%
for any $\phi \in C_{0}^{1}({\mathbb{R}}^{2})$. \ Since $0\leqslant \frac{1}{%
\delta }\int_{0}^{\delta }f(t,\cdot ,\cdot )\,dt\leqslant 1$ for $\delta >0,$
we conclude that $\mathbf{F}\cdot (1,0,0)|_{t=0}=f(0,\cdot ,\cdot )$ is a
bounded function in $L^{\infty }({\mathbb{R}}^{2}),$ which defines the trace 
$f^{0}$ at $t=0$\ of the function $f,$ which is defined as the limit 
\begin{equation*}
f^{0}=\lim_{\delta \rightarrow 0^{+}}\frac{1}{\delta }\int_{0}^{\delta
}f(t,\cdot ,\mathbf{\cdot })\,dt\text{\ \ a.e. on }{\mathbb{R}}^{2}.
\end{equation*}%
Obviously we have%
\begin{equation}
0\leqslant f^{0}\leqslant 1\quad \text{ \ \ \ \ a.e. \ in \ \ }{\mathbb{R}}%
^{2}.  \label{1}
\end{equation}%
Now let us choose $\psi =O_{\delta }(t)\,\phi (x,c)$ for $\phi \in C_{0}^{1}(%
{\mathbb{R}}^{2})$\ in \eqref{2.6}, where $O_{\delta }=1-1_{\delta }$ with 
\begin{equation*}
1_{\delta }(t)=%
\begin{cases}
1,\quad \text{if }t\geqslant \delta , \\ 
t/\delta ,\quad \text{if }0<t<\delta .\text{ }%
\end{cases}%
\end{equation*}%
Let us note that such function $\psi $\ can be taken in \eqref{2.6} using
the standard regularization procedure. Taking $\delta \rightarrow 0$ in the
obtained equality, we derive%
\begin{equation*}
\int_{{\mathbb{R}}^{2}}f^{0}(x,c)\,\phi (x,c)\,dxdc=\int_{{\mathbb{R}}^{2}}%
\mathrm{sign}_{+}\left( z_{0}(x)-c\right) \,\phi (x,c)\,dxdc\geqslant 0,
\end{equation*}%
that is 
\begin{equation}
f^{0}(x,s)=\mathrm{sign}_{+}\left( z_{0}(x)-c\right) \text{ \ \ \ \ \ \ for
a.e. }(x,c)\in {\mathbb{R}}^{2}.  \label{18_2}
\end{equation}

Applying the same reasoning to the vetorial function 
\begin{equation*}
\mathbf{F}(t,x,c)=\left( 1-f,\ (1-f)\left( c+\mathcal{W}\right) ,\ -\left(
(1-f)\left( \partial _{x}\mathrm{p}+\mathcal{W}\partial _{x}\mathcal{W}%
+c\partial _{x}\mathcal{W}\right) +\mu \right) \right) 
\end{equation*}%
and now choosing $\psi =O_{\delta }(t)\,\phi (x,c)$ for $\phi \in C_{0}^{1}({%
\mathbb{R}}^{2})$\ in \eqref{2.7}, we deduce 
\begin{equation}
(1-f^{0})(x,s)=\mathrm{sign}_{-}\left( z_{0}(x)-c\right) \text{ \ \ \ \ \ \
for a.e. }(x,c)\in {\mathbb{R}}^{2}.  \label{18_1}
\end{equation}

The identities \eqref{18_2} and \eqref{18_1} complete the claimed result %
\eqref{2.13}.

\hfill $\blacksquare $

\bigskip

\bigskip

\begin{lemma}
\label{lemma33} For $P-$a.e. $\omega\in\Omega$, there exists a function $%
v(t,x)=v(t,x, \omega),$ such that for $P-$a.e. $\omega\in\Omega$ 
\begin{equation*}
f(t,x,c)=\mathrm{sign}_{+}(v(t,x)-c)\text{\ \ \ \ \ \ for a.e. \ }(t,x,c)\in 
\mathbb{R}_{T}\times \mathbb{R}.
\end{equation*}
\end{lemma}

\noindent \textbf{Proof.} Let us introduce 
\begin{equation*}
\psi (t,x,c)=O_{\delta }(t)\varphi _{\varepsilon }(x)\varphi _{\varepsilon
}(c)\text{ \ \ \ \ \ \ for }(t,x,c)\in {\mathbb{R}}_{T}\times {\mathbb{R}}
\end{equation*}%
with $O_{\delta }(t)=(1{_{\delta }}(t)-1{_{\delta }}(t-t_{0}+\delta ))$ \
for\ $t_{0}\in (2\delta ,T)$ and $\varphi _{\varepsilon }(s)=1{_{\varepsilon
}}(s+\varepsilon ^{-1})-1{_{\varepsilon }}(s-\varepsilon ^{-1})$ \ for $s\in 
{\mathbb{R}}.$ Therefore choosing this $\psi $ as a test function in the
respective integral form of \eqref{RRR} and taking the limit on $\varepsilon
\rightarrow 0,$ with the help of (\ref{equ1}) and (\ref{equ2}), we get the
inequality for $P-$a.e. $\omega\in\Omega$%
\begin{equation}
\frac{1}{\delta }\int_{t_{0}-\delta }^{t_{0}}\int_{{\mathbb{R}}^{2}}F\
dtdxdc\leqslant \frac{1}{\delta }\int_{0}^{\delta }\int_{{\mathbb{R}}^{2}}F\
dtdxdc=A^{\delta }.  \label{ddF}
\end{equation}%
Since for $P-$a.e. $\omega\in\Omega$ 
\begin{equation*}
-\frac{1}{\delta }\int_{0}^{\delta }f^{2}(t)dt\leqslant -\left( \frac{1}{%
\delta }\int_{0}^{\delta }f(t)\ dt\right) ^{2},
\end{equation*}%
we get that $P-$a.e. $\omega\in\Omega$ 
\begin{equation*}
A^{\delta }\leqslant \int_{{\mathbb{R}}^{2}}\left[ \left( \frac{1}{\delta }%
\int_{0}^{\delta }f(t)\ dt\right) -\left( \frac{1}{\delta }\int_{0}^{\delta
}f(t)\ dt\right) ^{2}\right] \ dxdc.
\end{equation*}%
By $0\leqslant \frac{1}{\delta }\int_{0}^{\delta }f(t)\ dt\leqslant 1$ and
Lemma \ref{lemma2} the limit transition on $\delta \rightarrow 0$ and the
dominated convergence theorem implies%
\begin{equation*}
\limsup_{\delta \rightarrow 0^{+}}A^{\delta }\leqslant \int_{{\mathbb{R}}%
^{2}}(f^{0}-f{^{0}}^{2})\ dcdx=0
\end{equation*}%
and, consequently, the Lebesgue differentiation theorem and $%
F=f(1-f)\geqslant 0$ imply $\ \ \ $%
\begin{equation}
0\leqslant \int_{\mathbb{R}^{2}}F(t_{0},x,c)\ d{x}dc\leqslant 0\text{\ \ \ \
\ for a.e. }t_{0}\in \lbrack 0,T],
\end{equation}%
that is $F=0$ a.e. in $\mathbb{R}_{T}\times \mathbb{R}$ and $f$ \ can have
only the values $0$ and $1.$ Therefore, for $P-$a.e. $\omega\in\Omega$,
there exists a function $v=v(t,x, \omega),$ such that 
\begin{equation*}
f(t,x,c)=\mathrm{sign}_{+}(v(t,x)-c),
\end{equation*}%
since the function $f$ is monotone decreasing on $c$, by (\ref{eq2.88}%
).\hfill $\blacksquare $

\bigskip \bigskip

Therefore, for $P-$a.e. $\omega\in\Omega$, we have 
\begin{equation*}
\ |z_{\varepsilon }|_{+}=\int_{0}^{+\infty }f_{\varepsilon }(\mathbf{\cdot },%
\mathbf{\cdot },c)\ dc,\quad |z_{\varepsilon }|_{-}=\int_{-\infty
}^{0}(1-f_{\varepsilon }(\mathbf{\cdot },\mathbf{\cdot },c))\
dc\rightharpoonup |v|_{+},\quad |v|_{-}
\end{equation*}%
weakly -- $\ast $ in $L^{\infty }(0,T;L^{s}(\mathbb{R}))$ for any $s<q.$
This implies $z=|v|_{+}-|v|_{+}=v$ (see \eqref{eq2.88}).

The relations \eqref{limit0}$_{2},$\eqref{eq2.88}, \eqref{vvv} show that for 
$P-$a.e. $\omega\in\Omega$, it holds that 
\begin{equation*}
\Vert z_{\varepsilon }\Vert _{L^{s}(0,T;L^{s}(K))}\rightarrow \Vert z\Vert
_{L^{s}(0,T;L^{s}(K))}
\end{equation*}%
for any compact subset $K$ of $\mathbb{R}$. Adding the fact that for $P-$%
a.e. $\omega\in\Omega$ we have $z_{\varepsilon }\rightharpoonup z$ weakly$%
-\ast $ in ${L^{\infty }(0,T;L^{s}(\mathbb{R}))}$, for $P-$a.e. $%
\omega\in\Omega$ we obtain 
\begin{equation*}
z_{\varepsilon }\rightarrow z\quad \text{strongly in }L^{s}(0,T;L_{loc}^{s}(%
\mathbb{R})).
\end{equation*}%
Hence for $P-$a.e. $\omega\in\Omega$, we also have the strong convergence in 
${L^{1}(0,T;L_{loc}^{s}(\mathbb{R}))}$. Finally, using the H\"{o}lder
inequality, we can verify that the strong convergence of the sequence ${%
z_{\varepsilon }}$ in ${L^{1}(0,T;L_{loc}^{s}(\mathbb{R}))}$ and its
boundedness in $L^{\infty }(0,T;L^{s}(\mathbb{R}))$ imply 
\begin{equation*}
z_{\varepsilon }\rightarrow z\quad \text{strongly in }L^{r}(0,T;L_{loc}^{s}(%
\mathbb{R}))\quad \text{for}\quad 1\leq r<\infty .
\end{equation*}

\bigskip

The pair $z_{\varepsilon },$\ $\mathrm{p}_{\varepsilon }$\ satisfies %
\eqref{Sys}$_{2}$, and the equality obtained as the sum of \eqref{plus} and %
\eqref{minus} for $P-$a.e. $\omega \in \Omega $. The above strong
convergence and the convergences \eqref{limit0} imply that the limit
functions $z,\mathrm{p}$\ and $u=z+\mathcal{W}$ satisfy \eqref{1O}-%
\eqref{pdef} for $P-$a.e. $\omega \in \Omega $.

To finish the proof of Theorem \ref{the_1}, we apply Lemma 3.1 of \cite%
{Fla99} in order to construct a measurable solution of the system \eqref{1O}-%
\eqref{pdef} in the probabilistic space $(\Omega ,P)$. Namely, under the
assumptions of Theorem \ref{the_1}, we fix $1<r<\infty $, $s<q$, and
consider the Banach spaces 
\begin{align*}
Y& =\{y\in C([0,T];L^{2}(\mathbb{R})\cap L^{q}(\mathbb{R}))\;\text{with}%
\;\partial _{x}y\in L^{2}(0,T;L^{\infty }(\mathbb{R}))\}, \\
Z& =L^{r}(0,T;L^{2}(\mathbb{R})\cap L^{q}(\mathbb{R})\cap L^{s}(\mathbb{R}%
))\times L^{\infty }(0,T;L^{\infty }(\mathbb{R})),
\end{align*}%
endowed with the natural strong topologies. We have proved that for $P-$a.e. 
$\omega \in \Omega $ \ and $y=\mathcal{W}(\cdot ,\cdot ,\omega )\in Y$ \
there exists a pair $(u,\mathrm{p})\in Z,$ satisfying \eqref{1O}-\eqref{pdef}%
, that define the mapping 
\begin{equation*}
\Lambda :Y\rightarrow Z,\quad \Lambda (y)=(u,\mathrm{p}).
\end{equation*}%

This mapping is possibly multiple-valued, since in the limit transitions %
\eqref{limit0} we may select other converging sub-sequences and, as a
result, obtain some different limit functions $(u,\mathrm{p})$ with $u=z-%
\mathcal{W}$, satisfying the system \eqref{1O}-\eqref{pdef}.

Let us verify that the mapping $\Lambda $ has a closed graph. Taking a
sequence $(y_{n},\Lambda (y_{n}))$, we write the system \eqref{1O}-%
\eqref{pdef} with $(y,(u,\mathrm{p}))$ replaced by $(y_{n},(u_{n},%
\mathrm{p}_{n}))$ with $(u_{n},\mathrm{p}_{n})=\Lambda (y_{n}).$ \ Assuming
that 
\begin{equation*}
(y_{n},\Lambda (y_{n}))\rightarrow (y,\overline{\Lambda })\quad \text{%
strongly in }Y\times Z
\end{equation*}%
we easily pass to the limit in \eqref{1O}-\eqref{pdef}, as $n\rightarrow
\infty ,$ showing that $\overline{\Lambda }\equiv \Lambda (y)=(u,\mathrm{p})$%
\ also verifies this system with $\mathcal{W}$ replaced by $y.$

Therefore Lemma 3.1 in \cite{Fla99}
 implies that  there exists a measurable selection 
 $$
 U:Y\to Z, \qquad y\to (u,p).
 $$
 Since $\mathcal{W}:\Omega \to Y$ is measurable, 
 the solution of the   system (\ref%
{Hyp3}) is defined by the composition
  $
  (u,q)=U\circ \mathcal{W}:\Omega \to Z,
 $
 which is a measurable mapping.

\bigskip

\begin{remark}
We should notice that if, in addition,  the initial condition $u_0 $
belongs to $L^\infty(\mathbb{R})$, the deterministic problem has an unique solution, and consequently the pathwise uniqueness holds.
Futhermore, the solution is an adapted  stochastic process.
\end{remark}

\textbf{Acknowledgments}

\bigskip

Lynnyngs K. Arruda acknowledges support from FAPESP (Funda\c c\~ao de Amparo
\`a Pesquisa do Estado de S\~ao Paulo - Brazil), Grant 2017/23751-2.

The work of F. Cipriano was partially supported by the Funda\c{c}\~{a}o para
a Ci\^{e}ncia e a Tecnologia (Portuguese Foundation for Science and
Technology) through the project UID/MAT/00297/2019 (Centro de Matem\'{a}tica
e Aplica\c{c}\~{o}es).


%
\section*{Conflict of interest}
The authors declare that they have no conflict of interest.


%

\end{document}